\documentclass[reqno]{amsart}
\usepackage{setspace}
\usepackage{graphicx}
\usepackage{amsmath}
\usepackage{amsfonts}
\usepackage{amssymb}
\begin{document}

\centerline{\Huge \bf The Simultaneous Lattice Packing-Covering}
\centerline{\Huge \bf Constant of Octahedra}

\bigskip
\centerline{\large Yiming Li and Chuanming Zong}

{\large
\vspace{0.8cm}
\centerline{\begin{minipage}{12.8cm}
{\bf Abstract.} This paper proves that the simultaneous lattice packing-covering constant of an octahedron is $7/6$. In other words, $7/6$ is the smallest positive number $r$ such that for every octahedron $O$ centered at the origin there is a lattice $\Lambda$ such that $O+\Lambda$ is a packing in $\mathbb{E}^3$ and $rO+\Lambda$ is a covering of $\mathbb{E}^3$.
\end{minipage}}



\vspace{0.8cm}
\noindent
{\LARGE\bf 1. Introduction}

\bigskip
In 1950, C. A. Rogers \cite{Ro01} introduced and studied two simultaneous packing-covering constants $\gamma(C)$ and $\gamma^{*}(C)$ for an $n$-dimensional centrally symmetric convex body $C$ centered at the origin of $\mathbb{E}^n$:
$$\gamma(C)=\min_X\{r:\ rC+X\ {\rm is \ a \ covering \ of} \ {\mathbb E^n}\},$$
where $X$ is an arbitrary discrete point set such that $C+X$ is a packing in ${\mathbb E}^n$;
$$\gamma^{*}(C)=\min_\Lambda\{r:\ rC+\Lambda\ {\rm is \ a \ covering \ of} \ {\mathbb E^n}\},$$
where $\Lambda$ is a lattice such that $C+\Lambda$ is a packing in ${\mathbb E^n}$. By an inductive method, he proved that
$$\gamma (C)\le \gamma^{*}(C)\leq 3$$
holds for all $n$-dimensional centrally symmetric convex bodies. In 1972, via mean value techniques developed by C. A. Rogers \cite{Ro02} and C. L. Siegel \cite{Si01}, G. L. Butler \cite{Bu01} proved that
$$\gamma^{*}(C)\leq 2+o(1)$$
holds for all $n$-dimensional centrally symmetric convex bodies.

In 1970s, L. Fejes T\'{o}th \cite{Fe01, Fe02} introduced and investigated two deep hole constants $\rho(C)$ and $\rho^{*}(C)$ for an $n$-dimensional centrally symmetric convex body $C$ centered at the origin of $\mathbb{E}^n$, where $\rho(C)$ is the largest positive number $r$ such that one can put a translate of $rC$ into every translative packing $C+X$, and $\rho^{*}(C)$ is the largest positive number $r^{*}$ such that one can put a translate of $r^{*}C$ into every lattice packing $C+\Lambda$. Clearly, we have
$$\gamma(C)=\rho(C)+1$$
and
$$\gamma^{*}(C)=\rho^{*}(C)+1.$$

Let $B^{n}$ denote the $n$-dimensional unit ball centered at the origin. Like the packing density problem and the covering density problem, to determine the values of $\gamma(B^n)$ and $\gamma^{*}(B^n)$ is important and interesting. The known exact results are listed in the following table:

\medskip
\centerline{
\renewcommand\arraystretch{1.5}
\begin{tabular}{|c|c|c|c|c|c|}
\hline $n$ & $2$ & $3$ & $4$ & $5$ \\
\hline $\gamma^*(B^n)$ & $2/\sqrt3$ & $\sqrt{5/3}$ & $\sqrt{2\sqrt3}(\sqrt3-1)$ & $\sqrt{3/2+\sqrt{13}/6}$\\
\hline $\gamma(B^n)$ & $2/\sqrt3$ & $\sqrt{5/3}$ & ?? & ??\\
\hline Authors & Trivial & B\"{o}r\"{o}czky \cite{Bo01} & Horv\'{a}th \cite{Ho01} & Horv\'{a}th \cite{Ho02}\\
\hline
\end{tabular}}

\bigskip
In 1995, for the ball case, Rogers' reductive method was modified and his upper bound was improved by M. Henk \cite{He01} to
$$\gamma^*(B^n)\le \sqrt{21}/2= 2.29128\cdots.$$
Clearly, this upper bound is not as good as Butler's upper bound. However, Rogers' approach has application in lattice cryptography (see Micciancio \cite{Mi04}). On one hand, $\gamma^*(B^n)$ is a bridge connecting the shortest vector problem (SVP) and the closest vector problem (CVP), both fundamental in lattice cryptography. On the other hand, the reductive argument can lead to an algorithm.

In the plane, C. Zong \cite{Zo01, Zo03} proved that
$$\gamma(C)=\gamma^{*}(C)\leq 2(2-\sqrt2)$$
holds for all centrally symmetric convex domains and the second equality holds if and only if $C$ is an affinely regular octagon.
It is remarkable and interesting to notice that the maximum was not attained by circular discs! In ${\mathbb E}^3$,  C. Zong \cite{Zo02} proved that
$$\gamma^{*}(C)\leq 1.75$$
holds for all centrally symmetric convex bodies. It is interesting to compare with
$$\gamma (B^3)=\gamma^*(B^3)=\sqrt{5/3}=1.29099\cdots.$$

Let $O$ denote the regular octahedron with vertices $(1,0,0)$, $(0,1,0)$, $(0,0,1)$, $(-1,0,0)$, $(0,-1,0)$ and $(0,0,-1)$. In 1904, Minkowski \cite{Mi01} proved that: the lattice $\Lambda$ generated by ${\bf a}_1=(-\frac{2}{3},1,\frac{1}{3})$, ${\bf a}_2=(\frac{1}{3},-\frac{2}{3},1)$, ${\bf a}_3=(1,\frac{1}{3},-\frac{2}{3})$ give the optimal lattice packing density $18/19$. In fact, one can easily see that, Minkowski's result shows $\gamma^*(O)\leq7/6$ and all the uncovered spaces are regular tetrahedra, see Figure 1.

\begin{figure}[h]
\centering
\includegraphics[height=5.0cm]{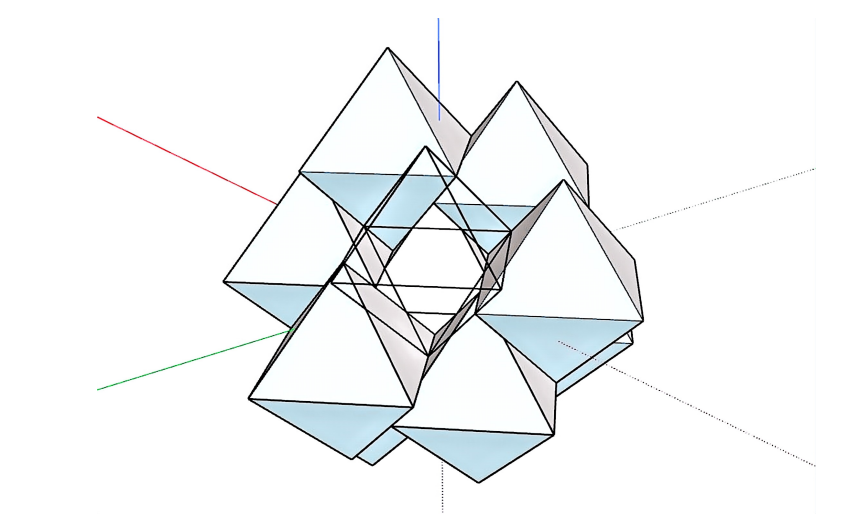}
\caption{The optimal lattice packing configuration for $O$.}
\label{fig1}
\end{figure}

By studying these tetrahedral holes and all its variations, this article will prove the following theorem:

\medskip
\noindent
{\bf Theorem 1.}
$$\gamma^{*}(O)=7/6.$$

\medskip
In section 4, we also investigate the simultaneous packing-covering constants of some other polytopes. Based on those examples, we propose the following problem:

\medskip
\noindent
{\bf Problem 1.} In ${\mathbb E}^3$, does
$$\gamma (C)=\gamma^{*}(C)\le \sqrt{5/3}$$
hold for every centrally symmetric convex body $C$?

\vspace{0.8cm}
\noindent
{\LARGE\bf 2. Technical Lemmas}

\bigskip
Let $T$ denote the regular tetrahedron with vertices (0,0,1), $({\frac{1}{3}},0,{\frac{2}{3}})$, $(0,{\frac{1}{3}},{\frac{2}{3}})$ and $({\frac{1}{3}},{\frac{1}{3}},1).$ Denote the distance function defined by $O$ as $||{\bf x}||_1$ (also known as $L_1$ measure), denote the distance between ${\bf x}$ and ${\bf y}$ in $L_1$-space as:
$$||{\bf x}, {\bf y}||_1=||{\bf x}-{\bf y}||_1 \ .$$

As usual, for two vector sets $A$ and $B$, define the minus set and Minkowski sum by
$$A\setminus B=\{{\bf x}:\ {\bf x}\in A\ {\rm and}\ {\bf x}\notin B\}$$
and
$$A+B=\{{\bf x}+{\bf y}:\ {\bf x}\in A\ {\rm and}\ {\bf y}\in B\},$$
respectively. We say a vector set $A$ is positively homothetic to $B$ with factor $r$, if $A=rB+{\bf x}$, $r>0$. We use ${\rm int}K$, ${\rm rint}K$, ${\rm cl}K$ and ${\rm conv}K$ to denote the interior of $K$, relative interior of $K$, the closure set of $K$ and the convex hull of $K$ as usual, use $\overline{{\bf x}{\bf y}}$ denote the segment with the vertices ${\bf x}$ and ${\bf y}$. We say $\gamma(C,X)=r$, if $C+X$ is a packing and $r$ is the minimum positive value which satisfy $rC+X$ is a covering.

Denote ${\bf e}_1=(1,0,0)$, ${\bf e}_2=(0,1,0)$, ${\bf e}_3=(0,0,1)$. Denote the faces of $O$:
\begin{align}
&{\rm conv}\{{\bf e}_1,{\bf e}_2,{\bf e}_3\},\ {\rm conv}\{{\bf e}_1,-{\bf e}_2,{\bf e}_3\},\ {\rm conv}\{-{\bf e}_1,-{\bf e}_2,{\bf e}_3\},\ {\rm conv}\{-{\bf e}_1,{\bf e}_2,{\bf e}_3\},\notag\\
&{\rm conv}\{{\bf e}_1,{\bf e}_2,-{\bf e}_3\},\ {\rm conv}\{{\bf e}_1,-{\bf e}_2,-{\bf e}_3\},\ {\rm conv}\{-{\bf e}_1,-{\bf e}_2,-{\bf e}_3\},\ {\rm conv}\{-{\bf e}_1,{\bf e}_2,-{\bf e}_3\}\notag
\end{align}
by $F^{(1)},F^{(2)},F^{(3)},F^{(4)},F^{(1')},F^{(2')},F^{(3')}$ and $F^{(4')}$, respectively.

\medskip
\noindent
{\bf Observation.} It is well known that, combining the regular tetrahedron with vertices ${\bf e}_1$, ${\bf e}_2,$ ${\bf e}_3$, ${\bf e}_1+{\bf e}_2+{\bf e}_3$, the regular tetrahedron with vertices $-{\bf e}_1$, $-{\bf e}_2,$ $-{\bf e}_3$, $-{\bf e}_1-{\bf e}_2-{\bf e}_3$ and $O$, we obtain a parallelepiped. Therefore, the sum of dihedral angle of a regular octahedron and a regular tetrahedron is $\pi$.

On the other hand, since a regular octahedron is defined by four pairs of parallel faces, we can observe that: if the intersections of each two of four regular octahedra $O$, $O+{\bf x}_1$, $O+{\bf x}_2$, $O+{\bf x}_3$ are two-dimensional, then the hole surrounded by them is a regular tetrahedron.

\medskip
It is natural to prove the following conclusion:

\medskip
\noindent
{\bf Lemma 1.} Suppose $F^{(1)}\supset O\cap (T+{\bf a}_1)\neq \emptyset$, then the center of gravity ${\bf g}$ of $T+{\bf a}_1$ satisfies $||{\bf g}||_1 \geq {\frac {7}{6}}$.

\begin{proof}
Define
$$X=\{{\bf g}:\ {\bf g}\ {\rm is\ the\ center\ of\ gravity\ of}\ T+{\bf a}_1,\ {\rm which\ satisfy}\ F^{(1)}\supset O\cap (T+{\bf a}_1)\neq \emptyset \}.$$
It is easy to see that
$$X\subset \{{(x,y,z):\ x+y+z=m}\},\ {\rm for\ a\ constant}\ m.$$

The center of gravity of $T$ is $(\frac{1}{6}, \frac{1}{6}, \frac{5}{6})$, therefore $m=\frac{7}{6}$.
Since $||{\bf o}, {\bf x}||_1\geq \frac{7}{6}$ holds for all ${\bf x}\in \big\{{(x,y,z):\ x+y+z=\frac{7}{6}}\big\}$, Lemma 1 is proved.
\end{proof}

\medskip
\noindent
{\bf Corollary 1.} For a packing $O+X$, suppose there is a regular tetrahedron $T_1$ satisfying the following two conditions:\\
$(1)$: ${\rm int}T_1\cap (O+X)= \emptyset,$\\
$(2)$: $T_1$ is positively homothetic to $T$ with dilation factor $r\geq1$,
\\
then we have $\gamma(O,X)\geq \frac{7}{6}$.

\begin{proof}
Denote the center of gravity of $T_1$ by ${\bf g}_1$. By Lemma 1, we have
$$||{\bf g}_1, {\bf x}||_1\geq\mbox{$\frac{7}{6}$},\ {\rm for\ all}\ {\bf x}\in X.$$
Therefore, ${\bf g}_1\notin {\rm int}(\frac{7}{6}O)+X,$ which means $\gamma(O,X)\geq\frac{7}{6}.$
\end{proof}

\medskip
Suppose ${\bf a}_0=(x_0,y_0,z_0)\in {\rm conv}\{(0,1,1),(0,\frac{4}{3},\frac{2}{3}),(\frac{1}{3},1,\frac{2}{3}),(\frac{1}{3},\frac{4}{3},1)\}$. In other words,
$$x_0+y_0+z_0\geq2,$$
$$-x_0+y_0+z_0\leq2,$$
$$x_0-y_0+z_0\leq0,$$
$$x_0+y_0-z_0\leq\mbox{$\frac{2}{3}$}. \eqno(1)$$

\medskip
For two convex bodies $C_1$ and $C_2$, we say $C_1$ is obstructed by $C_2$, if $C_2\cap {\rm int}(C_1)\neq\emptyset$. We will prove that, for a packing $O+X$ containing $O$ and $O+{\bf a}_0$, we have $\gamma(O,X)\geq\frac{7}{6}$. To this end, we show that no matter how to obstruct the unpacked place by the translative of $O$, there must exist a regular tetrahedron positively homothetic to $T$ with factor $r\geq1$, which is not obstructed by $O+X$.

\medskip
Define $T'={\rm conv}\{{\bf y}_1,{\bf y}_2,{\bf y}_3,{\bf y}_4\}$, where
\begin{eqnarray*}
& &{\bf y}_1=(1-\mbox{$\frac{-x_0+y_0+z_0}{2}$},y_0-1,1-\mbox{$\frac{x_0+y_0-z_0}{2}$}),\\
& &{\bf y}_2=(1-\mbox{$\frac{-x_0+y_0+z_0}{2}$},\mbox{$\frac{-x_0+y_0+z_0}{2}$},0),\\
& &{\bf y}_3=(2-y_0,y_0-1,0),\\
& &{\bf y}_4=(2-y_0,\mbox{$\frac{-x_0+y_0+z_0}{2}$},1-\mbox{$\frac{x_0+y_0-z_0}{2}$}).
\end{eqnarray*}
In other words, $T'$ is a regular tetrahedron positively homothetic to $T$ with factor $r_0$, formed by all the points $(x,y,z)$ which satisfy:
\begin{eqnarray*}
& &x+y+z\geq1,\\
& &-x+y+z\leq-x_0+y_0+z_0-1,\\
& &x-y+z\leq3-2y_0,\\
& &x+y-z\leq1.
\end{eqnarray*}
$$\eqno(2)$$

Since $1-\frac{x_0+y_0-z_0}{2}\geq\frac{2}{3}$ by (1), we have
$$r_0\geq2. \eqno(3)$$

\noindent
It is easy to see that, $T'$ contact both $O$ and $O+{\bf a}_0$ at its boundary.

\medskip
\noindent
{\bf Lemma 2.} Suppose a regular octahedron $O+{\bf a}_1$ satisfies $O+{\bf a}_1\cap {\rm int}(O\cup O+{\bf a}_0)=\emptyset$. Then the following two statements are equivalent:\\
$(1)$: $O+{\bf a}_1\cap {\rm int}T'\neq \emptyset$.\\
$(2)$: ${\bf a}_1\in T''+(1,0,0)$, where $T''=T' \setminus ({\rm conv}\{{\bf y}_1,{\bf y}_3,{\bf y}_4\}\cup {\rm conv}\{{\bf y}_2,{\bf y}_3,{\bf y}_4\})$.
\begin{proof}  Define $Y=\{{\bf a}_1:\ (O+{\bf a}_1)\cap {\rm int}T'\neq \emptyset\}$, by routine computation, we have
\begin{eqnarray*}
Y={\rm int}\Big({\rm conv}\hspace{-0.5cm} &\big\{ \hspace{-0.5cm}&{\bf y}_1+(0,0,1),{\bf y}_4+(0,0,1),{\bf y}_1+(0,-1,0),{\bf y}_1+(-1,0,0),\\
& &{\bf y}_4+(1,0,0), {\bf y}_4+(0,1,0),{\bf y}_3+(0,1,0),{\bf y}_3+(-1,0,0),\\
& &{\bf y}_2+(1,0,0),{\bf y}_2+(0,-1,0),{\bf y}_2+(0,0,-1),{\bf y}_3+(0,0,-1)\big\}\Big).
\end{eqnarray*}

Define
$$Y'=Y\setminus\{{\bf a}_1:\ O+{\bf a}_1\cap {\rm int}(O\cup O+{\bf a}_0)\neq\emptyset\},$$
then Lemma 2 holds if and only if $Y'=T''+(1,0,0)$.

On the one hand, we have
$$Y'=Y\setminus({\rm int}(2O)\cup {\rm int}(2O+{\bf a}_0)).$$
Since ${\bf y}_1,{\bf y}_2,{\bf y}_3\in O$, by the convexity of $2O$, we have
\begin{eqnarray*}
{\rm conv}\hspace{-0.5cm} &\{ \hspace{-0.5cm}&{\bf y}_1+(0,0,1), {\bf y}_1+(-1,0,0), {\bf y}_1+(0,-1,0), {\bf y}_2+(-1,0,0), {\bf y}_2+(0,1,0), \\
& &{\bf y}_2+(0,0,-1), {\bf y}_3+(1,0,0), {\bf y}_3+(0,-1,0), {\bf y}_3+(0,0,-1)\}
\end{eqnarray*}
is a subset of $2O$, denote this convex hull by $Y_1$.
Also, by routine computation, we have
\begin{eqnarray*}
& &{\bf y}_1+(1,0,0), {\bf y}_1+(0,0,1), {\bf y}_2+(1,0,0),{\bf y}_2+(0,1,0),\\
& &{\bf y}_4+(1,0,0),{\bf y}_4+(0,1,0), {\bf y}_4+(0,0,1)\in 2O+{\bf a}_0.
\end{eqnarray*}
For example,
\begin{eqnarray*}
||{\bf y}_4+(0,0,1)-{\bf a}_0\hspace{-0.3cm} &||_1 \hspace{-0.3cm}&=|2-x_0-y_0|+|\mbox{$\frac{-x_0-y_0+z_0}{2}$}|+|2-\mbox{$\frac{x_0+y_0+z_0}{2}$}|\\
& &=2-x_0-y_0-\mbox{$\frac{-x_0-y_0+z_0}{2}$}+2-\mbox{$\frac{x_0+y_0+z_0}{2}$}\\
& &=4-x_0-y_0-z_0\\
& &\leq2,
\end{eqnarray*}
by (1). Therefore, the convex hull
\begin{eqnarray*}
{\rm conv}\hspace{-0.5cm} &\{ \hspace{-0.5cm}&{\bf y}_1+(1,0,0), {\bf y}_1+(0,0,1), {\bf y}_2+(1,0,0),{\bf y}_2+(0,1,0), \\
& &{\bf y}_4+(1,0,0),{\bf y}_4+(0,1,0), {\bf y}_4+(0,0,1)\}
\end{eqnarray*}
is a subset of $2O+{\bf a}_0$, denote this convex hull by $Y_2$.
It is easy to see that
$$Y\subset {\rm int}Y_1\cup {\rm int}Y_2\cup(T''+(1,0,0)).$$
Therefore, we have
$$Y'\subset T''+(1,0,0). \eqno(4)$$

On the other hand, for vector ${\bf w}=(w_1,w_2,w_3)\in T''$, by (1) and (2), we have
$$w_1,w_2,w_3\geq0,\ w_1+w_2+w_3\geq1.$$
Thus, we have
$$||{\bf w}+(1,0,0), {\bf o}||_1=w_1+1+w_2+w_3\geq 2,$$
which means
$$O+{\bf w}+(1,0,0)\cap {\rm int}O=\emptyset.$$
Since ${\bf a}_0=(x_0,y_0,z_0)$ is on the plane $-x+y+z=-x_0+y_0+z_0$, and ${\bf w}+(1,0,0)$ is in the half space
$$-x+y+z\leq -x_0+y_0+z_0-2,$$
by (2), we have
$$O+{\bf w}+(1,0,0)\cap {\rm int}(O+{\bf a}_0)=\emptyset.$$

We argue as follows:

\medskip
\noindent
{\bf 1.} ${\bf w}\in {\rm int}T''$.
Combining with ${\bf w}\in {\rm int}T'\cap(O+{\bf w}+(1,0,0))\neq \emptyset$ obviously, by the definition of $Y'$, we have
$${\bf w}+(1,0,0)\in Y'.$$
{\bf 2.} ${\bf w}\in ({\rm conv}\{{\bf y}_1,{\bf y}_2,{\bf y}_4\}\cup {\rm conv}\{{\bf y}_1,{\bf y}_2,{\bf y}_3\})\setminus (\overline{{\bf y}_1{\bf y}_3}\cup\overline{{\bf y}_1{\bf y}_4}\cup\overline{{\bf y}_2{\bf y}_3}\cup\overline{{\bf y}_2{\bf y}_4}).$
Obviously: ${\rm conv}\{{\bf w},{\bf w}+(1,-1,0),{\bf w}+(2,0,0),{\bf w}+(1,1,0)\}$, a cross section of $O+{\bf w}+(1,0,0)$ intersects ${\rm int}T'$, which means ${\bf w}+(1,0,0)\in Y'$.

\medskip
\noindent
Therefore, we have $Y'\supset T''+(1,0,0)$. 

Combining with (4), we have $Y'=T''+(1,0,0)$, Lemma 2 is proved.
\end{proof}

\medskip
\noindent
{\bf Lemma 3.} Suppose a regular octahedron $O+{\bf a}_1$ satisfies $O+{\bf a}_1\cap {\rm int}(O\cup O+{\bf a}_0)=\emptyset$ and $O+{\bf a}_1\cap {\rm int}T'\neq \emptyset.$ Then for an arbitrary regular octahedron $O+{\bf a}_2$, the following two conditions cannot both hold:\\
$(a)$: $O+{\bf a}_2\cap {\rm int}(O\cup (O+{\bf a}_0)\cup (O+{\bf a}_1))= \emptyset$,\\
$(b)$: $O+{\bf a}_2\cap {\rm int}T'\neq \emptyset.$
\begin{proof} By Lemma 2, we have ${\bf a}_1\in T''+(1,0,0)$. If (a) and (b) hold simultaneously, we have
$${\bf a}_2\in T''+(1,0,0),$$
by Lemma 2. But for arbitrary points ${\bf x}, {\bf y}\in T''$, we must have
$$||{\bf x},{\bf y}||_1<2,$$
by (1) and (2). Therefore $O+{\bf a}_2\cap {\rm int}(O+{\bf a}_1)\neq \emptyset,$ which contradicts $(a)$.

Lemma 3 is proved.
\end{proof}

\medskip
\noindent
{\bf Corollary 2.} For a packing $O+X$ containing $O$ and $O+{\bf a}_0$, $\gamma(O,X)\geq\frac{7}{6}$ holds.
\begin{proof}
If ${\rm int}T'\cap (O+X)=\emptyset$, since $T'$ is positively homothetic to $T$ with factor $r_0\geq2$, by (3), the condition of Corollary 1 is satisfied.

Otherwise, if there exist ${\bf a}_1\in X$ satisfied ${\rm int}T'\cap (O+{\bf a}_1)\neq\emptyset$, by Lemma 3, $O+{\bf a}_1$ is the only regular octahedron in $O+X$ which intersect ${\rm int}T'$. In this case, there exist a regular tetrahedron $T'''\subset T'\setminus(O+{\bf a}_1)$ is positively homothetic to $T$ with factor $\frac{r_0}{2}\geq1$ and ${\rm int}T''\cap(O+X)=\emptyset$, which satisfies the condition of Corollary 1.
\end{proof}

\vspace{0.8cm}
\noindent
{\LARGE\bf 3. Proof of the Theorem}

\bigskip
Define
\begin{align}
&T^{(1)}_{1}=T,\quad T^{(1)}_{2}=T+(\mbox{$0,\frac{1}{3},-\frac{1}{3}$}),\quad T^{(1)}_{3}=T+(\mbox{$\frac{1}{3},0,-\frac{1}{3}$}),\notag\\
& T^{(1)}_{4}=T+(\mbox{$0,\frac{2}{3},-\frac{2}{3}$}),\quad T^{(1)}_{5}=T+(\mbox{$\frac{1}{3},\frac{1}{3},-\frac{2}{3}$}),\quad T^{(1)}_{6}=T+(\mbox{$\frac{2}{3},0,-\frac{2}{3}$}),\notag
\end{align}
respectively. It is easy to see that $T_i^{(1)}\cap O\subset F^{(1)},\ i\in\{1,\ldots,6\}$.

Rotate $T^{(1)}_i$ relative to $z-$axis anticlockwise for $\frac{\pi}{2}$ degree, $\pi$ degree, $\frac{3}{2}\pi$ degree, denote them by $T^{(2)}_i$, $T^{(3)}_i$, $T^{(4)}_i$, respectively. Define
$$T^{(1')}_{i}=-T^{(3)}_{i},\quad T^{(2')}_{i}=-T^{(4)}_{i},\quad T^{(3')}_{i}=-T^{(1)}_{i},\quad T^{(4')}_{i}=-T^{(2)}_{i}.$$ Similarly, we have $T_i^{(k)}\cap O\subset F^{(k)},\ i\in\{1,\ldots,6\},\ k\in\{1,2,3,4,1',2',3',4'\}$.

To generalize, we define $T^{(k+1)}_i$ by rotate $T^{(k)}_i$ relative to $z-$axis anticlockwise for $\frac{\pi}{2}$ degree. By the rotation, if $k_1\equiv k_2\ ({\rm mod}\ 4)$, then $T^{(k_1)}_i=T^{(k_2)}_i$.

In the centrally symmetric condition, we suppose $1', 2', 3', 4'$ is equivalent to $3, 4, 1, 2$, respectively.

Now we consider whether or not $T^{(k)}_i$ can be obstructed by a packing $O+X$, for all $i$, $k$. To this end, we define
$$P^{(k)}_{i}=\{{\bf x}:\ O+{\bf x}\cap {\rm int}T^{(k)}_{i}\neq \emptyset\ {\rm and}\ O+{\bf x}\cap {\rm int}O=\emptyset\},\ {\rm for\ all}\ i,k.$$
In fact,
$$P^{(k)}_{i}={\rm int}(T^{(k)}_{i}+O)\setminus{\rm int}(2O),\ {\rm for\ all}\ i,k.$$

\medskip
\noindent
By routine computation, we have
\begin{eqnarray*}
P^{(1)}_{1}\hspace{-0.2cm} &= \hspace{-0.2cm}& \ \ {\rm int}\Big({\rm conv}\big\{(\mbox{$0,0,2),(0,\frac{4}{3},\frac{2}{3}),(\frac{4}{3},0,\frac{2}{3}),(\frac{1}{3},\frac{1}{3},2),(\frac{1}{3},\frac{4}{3},1),(\frac{4}{3},\frac{1}{3},1$})\big\}\Big)\\
\hspace{-0.8cm} &\cup \hspace{-0.8cm}& \ \ {\rm rint}\Big({\rm conv}\big\{(\mbox{$0,0,2),(0,\frac{4}{3},\frac{2}{3}),(\frac{4}{3},0,\frac{2}{3}$})\big\}\Big),\\
P^{(1)}_{2}\hspace{-0.2cm} &= \hspace{-0.2cm}&P^{(1)}_{1}+(\mbox{$0,\frac{1}{3},-\frac{1}{3}$}),\quad P^{(1)}_{3}=P^{(1)}_{1}+(\mbox{$\frac{1}{3},0,-\frac{1}{3}$}),\\
P^{(1)}_{4}\hspace{-0.2cm} &= \hspace{-0.2cm}&P^{(1)}_{1}+(\mbox{$0,\frac{2}{3},-\frac{2}{3}$}),\quad P^{(1)}_{5}=P^{(1)}_{1}+(\mbox{$\frac{1}{3},\frac{1}{3},-\frac{2}{3}$}),\quad P^{(1)}_{6}=P^{(1)}_{1}+(\mbox{$\frac{2}{3},0,-\frac{2}{3}$}).
\end{eqnarray*}

\medskip
Define $M^{(k)}=P^{(k)}_{1}\cup P^{(k)}_{2}\cup P^{(k)}_{3}\cup P^{(k)}_{4}\cup P^{(k)}_{5}\cup P^{(k)}_{6}$ for all $k$. It is easy to see that $M^{(k_1)}\cap M^{(k_2)}=\emptyset$ for $k_1\neq k_2$.

\medskip
\noindent
We dissect $M^{(k)}$ into nineteen pieces as follow:

\begin{eqnarray*}
Q^{(k)}_1\hspace{-0.2cm} &= \hspace{-0.2cm}&\{{\bf x}:\ {\bf x}\in P^{(k)}_{1},\ {\bf x}\notin (P^{(k)}_{2}\cup P^{(k)}_{3}\cup P^{(k)}_{4}\cup P^{(k)}_{5}\cup P^{(k)}_{6})\},\\
Q^{(k)}_2\hspace{-0.2cm} &= \hspace{-0.2cm}&\{{\bf x}:\ {\bf x}\in (P^{(k)}_{1}\cap P^{(k)}_{2}),\ {\bf x}\notin (P^{(k)}_{3}\cup P^{(k)}_{4}\cup P^{(k)}_{5}\cup P^{(k)}_{6})\},\\
Q^{(k)}_3\hspace{-0.2cm} &= \hspace{-0.2cm}&\{{\bf x}:\ {\bf x}\in (P^{(k)}_{1}\cap P^{(k)}_{3}),\ {\bf x}\notin (P^{(k)}_{2}\cup P^{(k)}_{4}\cup P^{(k)}_{5}\cup P^{(k)}_{6})\},\\
Q^{(k)}_4\hspace{-0.2cm} &= \hspace{-0.2cm}&\{{\bf x}:\ {\bf x}\in (P^{(k)}_{1}\cap P^{(k)}_{2}\cap P^{(k)}_{3}\cap P^{(k)}_{5}),\ {\bf x}\notin  (P^{(k)}_{4}\cup P^{(k)}_{6})\},\\
Q^{(k)}_5\hspace{-0.2cm} &= \hspace{-0.2cm}&\{{\bf x}:\ {\bf x}\in (P^{(k)}_{1}\cap P^{(k)}_{2}\cap P^{(k)}_{4}),\ {\bf x}\notin (P^{(k)}_{3}\cup P^{(k)}_{5}\cup P^{(k)}_{6})\},\\
Q^{(k)}_6\hspace{-0.2cm} &= \hspace{-0.2cm}&\{{\bf x}:\ {\bf x}\in (P^{(k)}_{1}\cap P^{(k)}_{2}\cap P^{(k)}_{3}\cap P^{(k)}_{4}\cap P^{(k)}_{5}),\ {\bf x}\notin P^{(k)}_{6}\},\\
Q^{(k)}_7\hspace{-0.2cm} &= \hspace{-0.2cm}&\{{\bf x}:\ {\bf x}\in (P^{(k)}_{1}\cap P^{(k)}_{3}\cap P^{(k)}_{6}),\ {\bf x}\notin (P^{(k)}_{2}\cup P^{(k)}_{4}\cup P^{(k)}_{5})\},\\
Q^{(k)}_8\hspace{-0.2cm} &= \hspace{-0.2cm}&\{{\bf x}:\ {\bf x}\in (P^{(k)}_{1}\cap P^{(k)}_{2}\cap P^{(k)}_{3}\cap P^{(k)}_{5}\cap P^{(k)}_{6}),\ {\bf x}\notin  P^{(k)}_{4}\},\\
Q^{(k)}_9\hspace{-0.2cm} &= \hspace{-0.2cm}&\{{\bf x}:\ {\bf x}\in P^{(k)}_{4},\ {\bf x}\notin  (P^{(k)}_{1}\cup P^{(k)}_{2}\cup P^{(k)}_{3}\cup P^{(k)}_{5}\cup P^{(k)}_{6})\},\\
Q^{(k)}_{10}\hspace{-0.2cm} &= \hspace{-0.2cm}&\{{\bf x}:\ {\bf x}\in (P^{(k)}_{2}\cap P^{(k)}_{4}),\ {\bf x}\notin (P^{(k)}_{1}\cup P^{(k)}_{3}\cup P^{(k)}_{5}\cup P^{(k)}_{6})\},\\
Q^{(k)}_{11}\hspace{-0.2cm} &= \hspace{-0.2cm}&\{{\bf x}:\ {\bf x}\in (P^{(k)}_{4}\cap P^{(k)}_{5}),\ {\bf x}\notin (P^{(k)}_{1}\cup P^{(k)}_{2}\cup P^{(k)}_{3}\cup P^{(k)}_{6})\},\\
Q^{(k)}_{12}\hspace{-0.2cm} &= \hspace{-0.2cm}&\{{\bf x}:\ {\bf x}\in (P^{(k)}_{2}\cap P^{(k)}_{3}\cap P^{(k)}_{4}\cap P^{(k)}_{5}),\ {\bf x}\notin (P^{(k)}_{1}\cup  P^{(k)}_{6})\},\\
Q^{(k)}_{13}\hspace{-0.2cm} &= \hspace{-0.2cm}&\{{\bf x}:\ {\bf x}\in (P^{(k)}_{4}\cap P^{(k)}_{5}\cap P^{(k)}_{6}),\ {\bf x}\notin (P^{(k)}_{1}\cup P^{(k)}_{2}\cup P^{(k)}_{3})\},\\
Q^{(k)}_{14}\hspace{-0.2cm} &= \hspace{-0.2cm}&\{{\bf x}:\ {\bf x}\in (P^{(k)}_{2}\cap P^{(k)}_{3}\cap P^{(k)}_{4}\cap P^{(k)}_{5}\cap P^{(k)}_{6}),\ {\bf x}\notin P^{(k)}_{1}\},\\
Q^{(k)}_{15}\hspace{-0.2cm} &= \hspace{-0.2cm}&\{{\bf x}:\ {\bf x}\in P^{(k)}_{6},\ {\bf x}\notin (P^{(k)}_{1}\cup P^{(k)}_{2}\cup P^{(k)}_{3}\cup P^{(k)}_{4}\cup P^{(k)}_{5})\},\\
Q^{(k)}_{16}\hspace{-0.2cm} &= \hspace{-0.2cm}&\{{\bf x}:\ {\bf x}\in (P^{(k)}_{3}\cap P^{(k)}_{6}),\ {\bf x}\notin (P^{(k)}_{1}\cup P^{(k)}_{2}\cup P^{(k)}_{4}\cup P^{(k)}_{5})\},\\
Q^{(k)}_{17}\hspace{-0.2cm} &= \hspace{-0.2cm}&\{{\bf x}:\ {\bf x}\in (P^{(k)}_{5}\cap P^{(k)}_{6}),\ {\bf x}\notin (P^{(k)}_{1}\cup P^{(k)}_{2}\cup P^{(k)}_{3}\cup P^{(k)}_{4})\},\\
Q^{(k)}_{18}\hspace{-0.2cm} &= \hspace{-0.2cm}&\{{\bf x}:\ {\bf x}\in (P^{(k)}_{2}\cap P^{(k)}_{3}\cap P^{(k)}_{5}\cap P^{(k)}_{6}),\ {\bf x}\notin (P^{(k)}_{1}\cup P^{(k)}_{4})\},\\
Q^{(k)}_{19}\hspace{-0.2cm} &= \hspace{-0.2cm}&\{{\bf x}:\ {\bf x}\in (P^{(k)}_{2}\cap P^{(k)}_{3}\cap P^{(k)}_{5}),\ {\bf x}\notin (P^{(k)}_{1}\cup P^{(k)}_{4}\cup  P^{(k)}_{6})\},
\end{eqnarray*}
for all $k$. It is easy to see that $Q^{(k)}_{i_1}\cap Q^{(k)}_{i_2}=\emptyset$ for $i_1\neq i_2$, see Figure 2.

\begin{figure}[h]
\centering
\includegraphics[height=8.5cm]{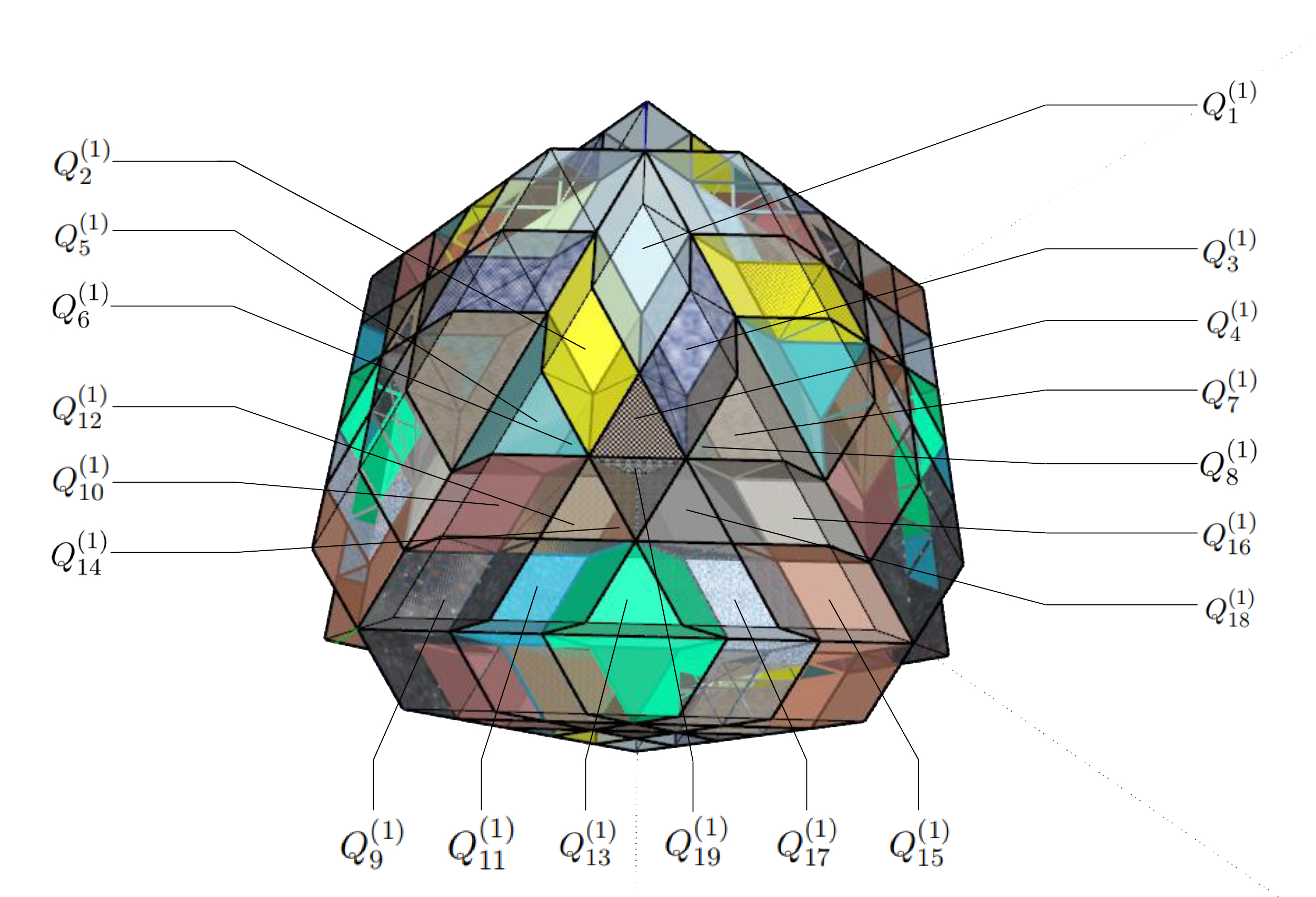}
\caption{Dissect $M^{(1)}$ into $Q^{(1)}_i$.}
\label{fig1}
\end{figure}

Furthermore, to show this decomposition more clearly, we give the following Figure 3.

\begin{figure}[h]
\centering
\includegraphics[height=6.5cm]{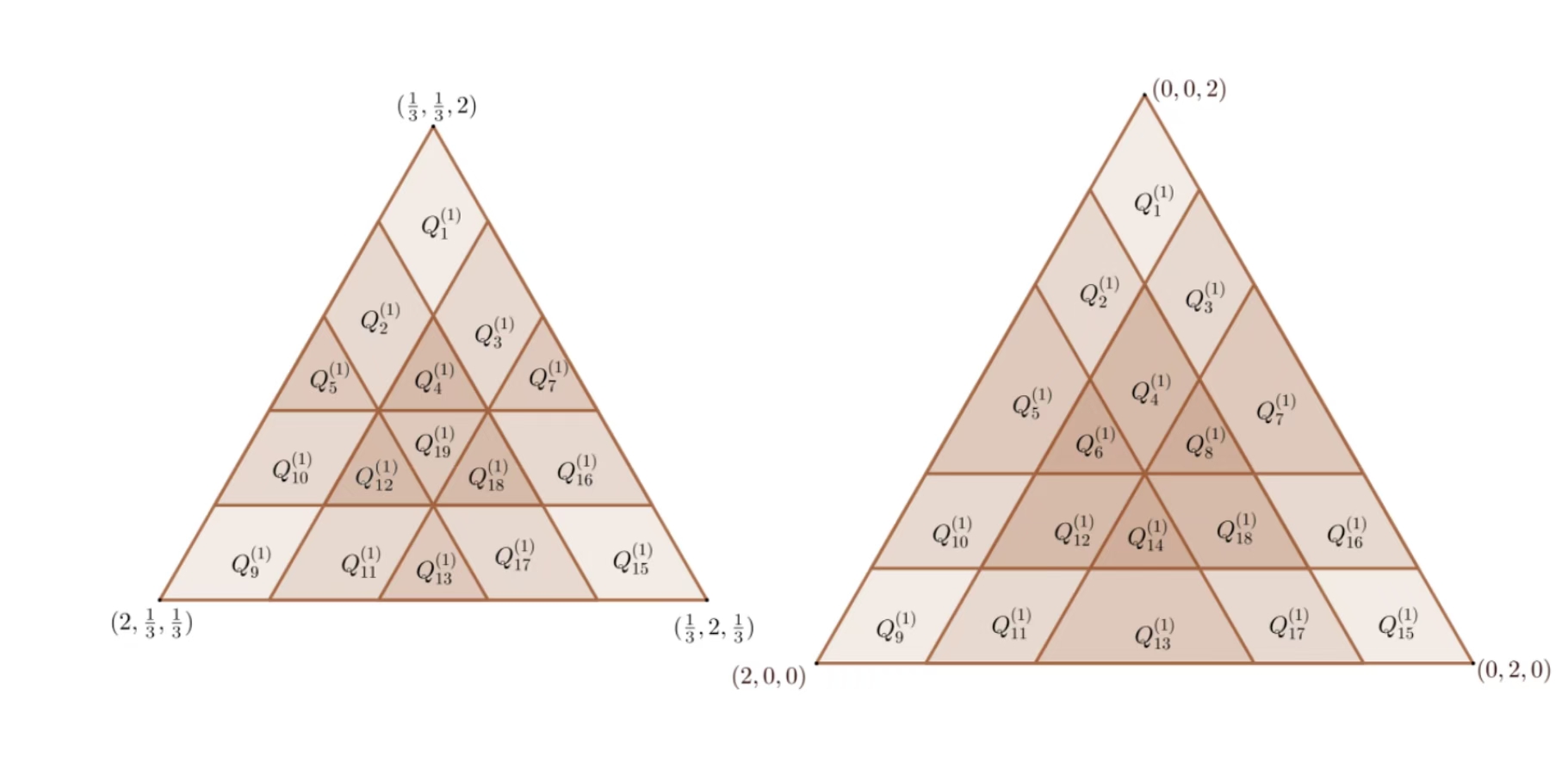}
\caption{The intersection of $M^{(1)}$ with the planes $\{x+y+z=8/3\}$ (left) and $\{x+y+z=2\}$ (right).}
\label{fig2}
\end{figure}

\noindent
We have:
\begin{eqnarray*}
{\rm cl}Q^{(1)}_1\hspace{-0.2cm} &= \hspace{-0.2cm}&{\rm conv}\big\{(\mbox{$0,0,2),(0,\frac{1}{3},\frac{5}{3}),(\frac{1}{3},\frac{1}{3},\frac{4}{3}),(\frac{1}{3},0,\frac{5}{3}),(\frac{1}{3},\frac{1}{3},2),(\frac{1}{3},\frac{2}{3},\frac{5}{3}),(\frac{2}{3},\frac{2}{3},\frac{4}{3}),(\frac{2}{3},\frac{1}{3},\frac{5}{3}$})\big\},\\
{\rm cl}Q^{(1)}_2\hspace{-0.2cm} &= \hspace{-0.2cm}&{\rm cl}Q^{(1)}_1+(\mbox{$0,\frac{1}{3},-\frac{1}{3}$}),\quad {\rm cl}Q^{(1)}_3={\rm cl}Q^{(1)}_1+(\mbox{$\frac{1}{3},0,-\frac{1}{3}$}),\\
{\rm cl}Q^{(1)}_4\hspace{-0.2cm} &= \hspace{-0.2cm}&{\rm conv}\big\{(\mbox{$\frac{1}{3},\frac{1}{3},\frac{4}{3}),(\frac{1}{3},\frac{2}{3},1),(\frac{2}{3},\frac{1}{3},1),(\frac{2}{3},\frac{2}{3},\frac{2}{3}),(\frac{2}{3},\frac{2}{3},\frac{4}{3}),(\frac{2}{3},1,1),(1,\frac{2}{3},1$})\big\},\\
{\rm cl}Q^{(1)}_5\hspace{-0.2cm} &= \hspace{-0.2cm}&{\rm conv}\big\{(\mbox{$0,\frac{2}{3},\frac{4}{3}),(\frac{1}{3},1,\frac{4}{3}),(\frac{1}{3},\frac{2}{3},1),(\frac{2}{3},1,1),(0,\frac{4}{3},\frac{2}{3}),(\frac{1}{3},\frac{4}{3},1),(\frac{1}{3},1,\frac{2}{3}$})\big\},\\
{\rm cl}Q^{(1)}_6\hspace{-0.2cm} &= \hspace{-0.2cm}&{\rm conv}\big\{(\mbox{$\frac{1}{3},\frac{2}{3},1),(\frac{2}{3},1,1),(\frac{1}{3},1,\frac{2}{3}),(\frac{2}{3},\frac{2}{3},\frac{2}{3}$})\big\},\\
{\rm cl}Q^{(1)}_7\hspace{-0.2cm} &= \hspace{-0.2cm}&{\rm conv}\big\{(\mbox{$\frac{2}{3},0,\frac{4}{3}),(1,\frac{1}{3},\frac{4}{3}),(\frac{2}{3},\frac{1}{3},1),(1,\frac{2}{3},1),(\frac{4}{3},\frac{1}{3},1),(1,\frac{1}{3},\frac{2}{3}),(\frac{4}{3},0,\frac{2}{3}$})\big\},\\
{\rm cl}Q^{(1)}_8\hspace{-0.2cm} &= \hspace{-0.2cm}&{\rm conv}\big\{(\mbox{$\frac{2}{3},\frac{1}{3},1),(1,\frac{2}{3},1),(1,\frac{1}{3},\frac{2}{3}),(\frac{2}{3},\frac{2}{3},\frac{2}{3}$})\big\},\\
{\rm cl}Q^{(1)}_9\hspace{-0.2cm} &= \hspace{-0.2cm}&{\rm conv}\big\{(\mbox{$0,2,0),(\frac{1}{3},\frac{5}{3},0),(\frac{1}{3},2,\frac{1}{3}),(\frac{2}{3},\frac{5}{3},\frac{1}{3}),(0,\frac{5}{3},\frac{1}{3}),(\frac{1}{3},\frac{4}{3},\frac{1}{3}),(\frac{1}{3},\frac{5}{3},\frac{2}{3}),(\frac{2}{3},\frac{4}{3},\frac{2}{3}$})\big\},\\
{\rm cl}Q^{(1)}_{10}\hspace{-0.2cm} &= \hspace{-0.2cm}&{\rm cl}Q^{(1)}_9+(\mbox{$0,-\frac{1}{3},\frac{1}{3}$}),\quad {\rm cl}Q^{(1)}_{11}={\rm cl}Q^{(1)}_9+(\mbox{$\frac{1}{3},-\frac{1}{3},0$}),\\
{\rm cl}Q^{(1)}_{12}\hspace{-0.2cm} &= \hspace{-0.2cm}&{\rm conv}\big\{(\mbox{$\frac{1}{3},1,\frac{2}{3}),(\frac{2}{3},\frac{2}{3},\frac{2}{3}),(\frac{1}{3},\frac{4}{3},\frac{1}{3}),(\frac{2}{3},1,\frac{1}{3}),(\frac{2}{3},1,1),(\frac{2}{3},\frac{4}{3},\frac{2}{3}),(1,1,\frac{2}{3}$})\big\},\\
{\rm cl}Q^{(1)}_{13}\hspace{-0.2cm} &= \hspace{-0.2cm}&{\rm conv}\big\{(\mbox{$\frac{2}{3},1,\frac{1}{3}),(1,\frac{2}{3},\frac{1}{3}),(\frac{2}{3},\frac{4}{3},0),(\frac{4}{3},\frac{2}{3},0),(1,1,\frac{2}{3}),(1,\frac{4}{3},\frac{1}{3}),(\frac{4}{3},1,\frac{1}{3}$})\big\},\\
{\rm cl}Q^{(1)}_{14}\hspace{-0.2cm} &= \hspace{-0.2cm}&{\rm conv}\big\{(\mbox{$\frac{2}{3},\frac{2}{3},\frac{2}{3}),(\frac{2}{3},1,\frac{1}{3}),(1,\frac{2}{3},\frac{1}{3}),(1,1,\frac{2}{3}$})\big\},\\
{\rm cl}Q^{(1)}_{15}\hspace{-0.2cm} &= \hspace{-0.2cm}&{\rm conv}\big\{(\mbox{$\frac{4}{3},\frac{1}{3},\frac{1}{3}),(\frac{5}{3},0,\frac{1}{3}),(\frac{5}{3},\frac{1}{3},0),(2,0,0),(\frac{4}{3},\frac{2}{3},\frac{2}{3}),(\frac{5}{3},\frac{1}{3},\frac{2}{3}),(\frac{5}{3},\frac{2}{3},\frac{1}{3}),(2,\frac{1}{3},\frac{1}{3}$})\big\},\\
{\rm cl}Q^{(1)}_{16}\hspace{-0.2cm} &= \hspace{-0.2cm}&{\rm cl}Q^{(1)}_{15}+(\mbox{$-\frac{1}{3},0,\frac{1}{3}$}),\quad {\rm cl}Q^{(1)}_{17}={\rm cl}Q^{(1)}_{15}+(\mbox{$-\frac{1}{3},\frac{1}{3},0$}),\\
{\rm cl}Q^{(1)}_{18}\hspace{-0.2cm} &= \hspace{-0.2cm}&{\rm conv}\big\{(\mbox{$\frac{2}{3},\frac{2}{3},\frac{2}{3}),(1,\frac{1}{3},\frac{2}{3}),(1,\frac{2}{3},\frac{1}{3}),(\frac{4}{3},\frac{1}{3},\frac{1}{3}),(1,\frac{2}{3},1),(1,1,\frac{2}{3}),(\frac{4}{3},\frac{2}{3},\frac{2}{3}$})\big\},\\
{\rm cl}Q^{(1)}_{19}\hspace{-0.2cm} &= \hspace{-0.2cm}&{\rm conv}\big\{(\mbox{$\frac{2}{3},\frac{2}{3},\frac{2}{3}),(\frac{2}{3},1,1),(1,\frac{2}{3},1),(1,1,\frac{2}{3}$})\big\}.
\end{eqnarray*}

\medskip
\noindent
{\bf Theorem 1.} $$\gamma^{*}(O)=7/6.$$
\begin{proof} Suppose $O+\Lambda$ is a lattice packing and $\gamma(O, \Lambda)<  \frac{7}{6}$. By Corollary 1, we have
$$(O+\Lambda) \cap {\rm int}T^{(k)}_i \neq \emptyset,\ {\rm for\ all}\ i,k.$$
By the definition of $P^{(k)}_i$, we have
$$\Lambda \cap P^{(k)}_i \neq \emptyset,\ {\rm for\ all}\ i,k.\eqno(5)$$

Obviously, we have
\begin{align}
&P^{(k)}_1=Q^{(k)}_1\cup Q^{(k)}_2\cup Q^{(k)}_3\cup Q^{(k)}_4\cup Q^{(k)}_5\cup Q^{(k)}_6\cup Q^{(k)}_7\cup Q^{(k)}_8,\tag{6}\\
&P^{(k)}_2=Q^{(k)}_2\cup Q^{(k)}_4\cup Q^{(k)}_5\cup Q^{(k)}_6\cup Q^{(k)}_8\cup Q^{(k)}_{10}\cup Q^{(k)}_{12}\cup Q^{(k)}_{14}\cup Q^{(k)}_{18}\cup Q^{(k)}_{19},\tag{7}\\
&P^{(k)}_3=Q^{(k)}_3\cup Q^{(k)}_4\cup Q^{(k)}_6\cup Q^{(k)}_7\cup Q^{(k)}_8\cup Q^{(k)}_{12}\cup Q^{(k)}_{14}\cup Q^{(k)}_{16}\cup Q^{(k)}_{18}\cup Q^{(k)}_{19},\tag{8}\\
&P^{(k)}_4=Q^{(k)}_5\cup Q^{(k)}_6\cup Q^{(k)}_9\cup Q^{(k)}_{10}\cup Q^{(k)}_{11}\cup Q^{(k)}_{12}\cup Q^{(k)}_{13}\cup Q^{(k)}_{14},\tag{9}\\
&P^{(k)}_5=Q^{(k)}_4\cup Q^{(k)}_6\cup Q^{(k)}_8\cup Q^{(k)}_{11}\cup Q^{(k)}_{12}\cup Q^{(k)}_{13}\cup Q^{(k)}_{14}\cup Q^{(k)}_{17}\cup Q^{(k)}_{18}\cup Q^{(k)}_{19},\tag{10}\\
&P^{(k)}_6=Q^{(k)}_7\cup Q^{(k)}_8\cup Q^{(k)}_{13}\cup Q^{(k)}_{14}\cup Q^{(k)}_{15}\cup Q^{(k)}_{16}\cup Q^{(k)}_{17}\cup Q^{(k)}_{18}.\tag{11}
\end{align}
Since $O+\Lambda$ is a lattice packing, by dilating the lattice a little we can keep that $\gamma(O,\Lambda)<\frac{7}{6}$ still holds and
$$||{\bf x},{\bf y}||_1>2,\ {\rm for\ all}\ {\bf x}, {\bf y}\in\Lambda,\ {\bf x}\neq{\bf y}.$$
In other words,
$${\bf y}\notin(2O+{\bf x}),\ {\rm for\ all}\ {\bf x}, {\bf y}\in\Lambda,\ {\bf x}\neq{\bf y}.$$

Denote $\mathcal{Q}_1=\{1,9,15\}$, $\mathcal{Q}_2=\{2,3,10,11,16,17\}$, $\mathcal{Q}_3=\{5,7,13,19\}$, $\mathcal{Q}_4=\{4,12,18\}$, $\mathcal{Q}_5=\{6,8,14\}$, which means that if $m\in \mathcal{Q}_j$, then $Q^{(k)}_m$ is completely inside exactly $j$ of $P^{(k)}_i$. Denote $\Lambda^{(k)}=\{i_1,i_2,...,i_n\}$, if
$$\Lambda\cap Q^{(k)}_{i}\neq \emptyset,\ {\rm for}\ i\in\{i_1,i_2,...,i_n\}$$
and
$$\Lambda\cap Q^{(k)}_{i}= \emptyset,\ {\rm for}\ i\notin\{i_1,i_2,...,i_n\}.$$

\medskip
\noindent
{\bf 3.1. All the possible $\Lambda^{(k)}$ for a given $k$}

To enumerate all the possible $\Lambda^{(k)}$ which satisfies (5) for a given $k$, we list some restricting conditions as follow:

For an arbitrary point ${\bf x}\in Q^{(k)}_{2}$, we have
$$(2O+{\bf x})\supset \big(Q^{(k)}_{3}\cup Q^{(k)}_{4}\cup Q^{(k)}_{5}\cup Q^{(k)}_{6}\cup Q^{(k)}_{8}\big).$$
Therefore, if $2\in\Lambda^{(k)}$, we must have
$$\Lambda^{(k)} \cap \{3,4,5,6,8\}= \emptyset. \eqno(12)$$
By symmetry, we have:
\begin{align}
&{\rm if}\ 3\in\Lambda^{(k)},\ {\rm then}\ \Lambda^{(k)} \cap \{2,4,6,7,8\}= \emptyset, \tag{13}\\
&{\rm if}\ 10\in\Lambda^{(k)},\ {\rm then}\ \Lambda^{(k)} \cap \{5,6,11,12,14\}= \emptyset, \tag{14}\\
&{\rm if}\ 11\in\Lambda^{(k)},\ {\rm then}\ \Lambda^{(k)} \cap \{6,10,12,13,14\}= \emptyset, \tag{15}\\
&{\rm if}\ 16\in\Lambda^{(k)},\ {\rm then}\ \Lambda^{(k)} \cap \{7,8,14,17,18\}= \emptyset,\tag{16}\\
&{\rm if}\ 17\in\Lambda^{(k)},\ {\rm then}\ \Lambda^{(k)} \cap \{8,13,14,16,18\}= \emptyset. \tag{17}
\end{align}

For an arbitrary point ${\bf x}\in Q^{(k)}_{4}$, we have
$$(2O+{\bf x})\supset \big(Q^{(k)}_{2}\cup Q^{(k)}_{3}\cup Q^{(k)}_{5}\cup Q^{(k)}_{6}\cup Q^{(k)}_{7}\cup Q^{(k)}_{8}\cup Q^{(k)}_{12}\cup Q^{(k)}_{14}\cup Q^{(k)}_{18}\cup Q^{(k)}_{19}\big).$$
Therefore, if $4\in\Lambda^{(k)}$, we must have
$$\Lambda^{(k)} \cap \{2,3,5,6,7,8,12,14,18,19\}= \emptyset. \eqno(18)$$
By symmetry, we have:
$${\rm if}\ 12\in\Lambda^{(k)},\ {\rm then}\ \Lambda^{(k)} \cap \{4,5,6,8,10,11,13,14,18,19\}= \emptyset, \eqno(19)$$
$${\rm if}\ 18\in\Lambda^{(k)},\ {\rm then}\ \Lambda^{(k)} \cap \{4,6,7,8,12,13,14,16,17,19\}= \emptyset.\eqno(20)$$

For an arbitrary point ${\bf x}\in Q^{(k)}_{5}$, we have
$$(2O+{\bf x})\supset \big(Q^{(k)}_{2}\cup Q^{(k)}_{4}\cup Q^{(k)}_{6}\cup Q^{(k)}_{8}\cup Q^{(k)}_{10}\cup Q^{(k)}_{12}\cup Q^{(k)}_{14}\cup Q^{(k)}_{19}\big).$$
Therefore, if $5\in\Lambda^{(k)}$, we must have
$$\Lambda^{(k)} \cap \{2,4,6,8,10,12,14,19\}= \emptyset. \eqno(21)$$
By symmetry, we have:
\begin{align}
&{\rm if}\ 7\in\Lambda^{(k)},\ {\rm then}\ \Lambda^{(k)} \cap \{3,4,6,8,14,16,18,19\}= \emptyset, \tag{22}\\
&{\rm if}\ 13\in\Lambda^{(k)},\ {\rm then}\ \Lambda^{(k)} \cap \{6,8,11,12,14,17,18,19\}= \emptyset. \tag{23}
\end{align}

For an arbitrary point ${\bf x}\in Q^{(k)}_{6}$, we have
$$(2O+{\bf x})\supset (Q^{(k)}_{2}\cup Q^{(k)}_{3}\cup Q^{(k)}_{4}\cup Q^{(k)}_{5}\cup Q^{(k)}_{7}\cup Q^{(k)}_{8}\cup Q^{(k)}_{10}\cup Q^{(k)}_{11}\cup Q^{(k)}_{12}\cup Q^{(k)}_{13}\cup Q^{(k)}_{14}\cup Q^{(k)}_{18}\cup Q^{(k)}_{19}).$$
Therefore, if $6\in\Lambda^{(k)}$, we must have
$$\Lambda^{(k)} \cap \{2,3,4,5,7,8,10,11,12,13,14,18,19\}= \emptyset. \eqno(24)$$
By symmetry, we have:
\begin{align}
&{\rm if}\ 8\in\Lambda^{(k)},\ {\rm then}\ \Lambda^{(k)} \cap \{2,3,4,5,6,7,12,13,14,16,17,18,19\}= \emptyset, \tag{25}\\
&{\rm if}\ 14\in\Lambda^{(k)},\ {\rm then}\ \Lambda^{(k)} \cap \{4,5,6,7,8,10,11,12,13,16,17,18,19\}= \emptyset. \tag{26}
\end{align}

For an arbitrary point ${\bf x}\in Q^{(1)}_{3}$, we have
$$(2O+{\bf x})\supset {\rm conv}\big\{(\mbox{$\frac{1}{3},\frac{2}{3},1),(0,\frac{2}{3},\frac{4}{3}),(0,1,1),(\frac{1}{3},1,\frac{2}{3}),(\frac{1}{3},1,\frac{4}{3}),(\frac{1}{3},\frac{4}{3},1),(\frac{2}{3},1,1$})\big\}.$$
Combining with Corollary 2, we have:
$${\rm if}\ \Lambda\cap Q^{(1)}_{3}\neq \emptyset,\ {\rm then}\ \Lambda\cap Q^{(1)}_{5}= \emptyset.$$

For an arbitrary point
$${\bf x}\in {\rm conv}\big\{(\mbox{$1,\frac{1}{3},\frac{2}{3}),(\frac{2}{3},0,\frac{4}{3}),(1,0,1),(\frac{2}{3},\frac{1}{3},1),(1,\frac{1}{3},\frac{4}{3}),(\frac{4}{3},\frac{1}{3},1),(1,\frac{2}{3},1$})\big\},$$
we have
$$(2O+{\bf x})\supset {\rm conv}\big\{(\mbox{$\frac{1}{3},\frac{2}{3},1),(0,\frac{2}{3},\frac{4}{3}),(0,1,1),(\frac{1}{3},1,\frac{2}{3}),(\frac{1}{3},1,\frac{4}{3}),(\frac{1}{3},\frac{4}{3},1),(\frac{2}{3},1,1$})\big\}.$$
By the symmetry of $2O$ and $M^{(1)}$, combining with Corollary 2, we have
$${\rm if}\ \Lambda\cap Q^{(1)}_{7}\neq \emptyset,\ {\rm then}\ \Lambda\cap Q^{(1)}_{5}= \emptyset.$$
By symmetry, we have:
\begin{align}
&{\rm if}\ 5\in\Lambda^{(k)},\ {\rm then}\ \Lambda^{(k)}\cap \{3,7,11,13\}= \emptyset, \tag{27}\\
&{\rm if}\ 7\in\Lambda^{(k)},\ {\rm then}\ \Lambda^{(k)}\cap \{2,5,13,17\}= \emptyset,\tag{28}\\
&{\rm if}\ 13\in\Lambda^{(k)},\ {\rm then}\ \Lambda^{(k)}\cap \{5,7,10,16\}= \emptyset.\tag{29}
\end{align}

Without loss of generality, if for two different lattices $\Lambda_1$ and $\Lambda_2$ satisfy (5), we have $\Lambda^{(k)}_1\subset\Lambda^{(k)}_2$, then we only consider $\Lambda^{(k)}_1$ instead of both. Suppose
$$i_1\in\mathcal{Q}_{j_1},\ i_2\in\mathcal{Q}_{j_2}, ... ,\ i_n\in\mathcal{Q}_{j_n},$$
to satisfying (5), a necessary condition is
$$j_1+j_2+...+j_n\geq6.$$
Combining with conditions (5)-(29), we categorize all the possible $\Lambda^{(k)}$ for a given $k$ as follow:

\medskip
\noindent
{\bf Category 1.} $n=2$ and $i_1\in\mathcal{Q}_5$, $i_2\in\mathcal{Q}_1$.

For instance, let $i_1=6$. By (5), (11) and (24), we have $i_2=15$. Therefore,
$$\Lambda^{(k)}=\{14,1\}, \{6,15\}, \{8,9\},$$
by the symmetry of $2O$ and $M^{(k)}$.

\medskip
\noindent
{\bf Category 2.} $n=2$ and $i_1\in\mathcal{Q}_5$, $i_2\in\mathcal{Q}_2$.

For instance, let $i_1=6$. By (5), (11) and (24), we have $i_2=16\ {\rm or}\ 17$. Therefore,
$$\Lambda^{(k)}=\{14,2\}, \{14,3\}, \{8,10\}, \{8,11\}, \{6,16\}, \{6,17\},$$
by the symmetry of $2O$ and $M^{(k)}$.

\medskip
\noindent
{\bf Category 3.} $n=2$ and $i_1\in\mathcal{Q}_4$, $i_2\in\mathcal{Q}_3$.

For instance, let $i_1=4$. By (5), (9), (11) and (18), we have $i_2=13$. Therefore,
$$\Lambda^{(k)}=\{4,13\}, \{18,5\}, \{12,7\},$$
by the symmetry of $2O$ and $M^{(k)}$.

\medskip
\noindent
{\bf Category 4.} $n=3$ and $i_1\in\mathcal{Q}_4$, $i_2,i_3\in\mathcal{Q}_1$.

For instance, let $i_1=4$. By (5), (9), (11) and (18), we have $i_2=9, i_3=15$. Therefore,
$$\Lambda^{(k)}=\{1,9,18\}, \{1,12,15\}, \{4,9,15\},$$
by the symmetry of $2O$ and $M^{(k)}$.

\medskip
\noindent
{\bf Category 5.} $n=3$ and $i_1\in\mathcal{Q}_4$, $i_2\in\mathcal{Q}_2$ satisfies ${\rm cl}Q^{(k)}_{i_{j_1}}\cap {\rm cl}Q^{(k)}_{i_{j_2}}\neq\emptyset$.

For instance, let $i_1=4$, $i_2=10$. By (5), (11), (14) and (18), we have $i_3=15\ {\rm or}\ 16\ {\rm or}\ 17$. Therefore,
\begin{eqnarray*}
\Lambda^{(k)}\hspace{-0.2cm} &= \hspace{-0.2cm}&\{4,10,15\}, \{4,10,16\}, \{4,10,17\}, \{4,16,9\}, \{4,16,11\},\\
& &\{12,2,15\}, \{12,2,16\}, \{12,2,17\}, \{12,17,1\}, \{12,17,3\},\\
& &\{18,11,1\}, \{18,11,2\}, \{18,11,3\}, \{18,3,10\}, \{18,3,9\},
\end{eqnarray*}
by the symmetry of $2O$ and $M^{(k)}$.

\medskip
\noindent
{\bf Category 6.} $n=3$ and $i_1\in\mathcal{Q}_4$, $i_2\in\mathcal{Q}_2$, satisfies ${\rm cl}Q^{(k)}_{i_1}\cap \big({\rm cl}Q^{(k)}_{i_2}\cup {\rm cl}Q^{(k)}_{i_3}\big)=\emptyset$.

For instance, let $i_1=4$, $i_2=11$. By (5), (11), (15) and (18), we have $i_3=15\ {\rm or}\ 17$. Therefore,
\begin{eqnarray*}
\Lambda^{(k)}\hspace{-0.2cm} &= \hspace{-0.2cm}&\{4,11,15\}, \{4,11,17\}, \{4,17,9\},\\
& &\{12,3,15\}, \{12,3,16\}, \{12,16,1\},\\
& &\{18,10,2\}, \{18,10,1\}, \{18,2,9\},
\end{eqnarray*}
by the symmetry of $2O$ and $M^{(k)}$.

\medskip
\noindent
{\bf Category 7.} $n=3$ and $i_1,i_2,i_3\in\mathcal{Q}_2$. By (5)-(17), it is easy to deduce that
$$\Lambda^{(k)}=\{2,11,16\}, \{3,10,17\}.$$

\medskip
\noindent
{\bf Category 8.} $n=4$ and $i_1=19$, $i_2\in\mathcal{Q}_2$.

For instance, let $i_2=2$. By (5), (9), (11) and (12), we have $i_3=9\ {\rm or}\ 10\ {\rm or}\ 11$ and $i_4=15\ {\rm or}\ 16\ {\rm or}\ 17$. Therefore,
\begin{eqnarray*}
\Lambda^{(k)}\hspace{-0.2cm} &= \hspace{-0.2cm}&\{19,2,9,15\}, \{19,2,9,16\}, \{19,2,9,17\},  \{19,2,10,15\}, \{19,2,10,16\},\\
& &\{19,2,10,17\}, \{19,2,11,15\}, \{19,2,11,17\}, \{19,3,9,15\}, \{19,3,9,16\},\\
& &\{19,3,9,17\}, \{19,3,10,15\}, \{19,3,10,16\}, \{19,3,11,15\}, \{19,3,11,16\},\\
& &\{19,3,11,17\}, \{19,16,9,1\}, \{19,16,10,1\}, \{19,16,11,1\}, \{19,17,9,1\},\\
& &\{19,17,10,1\}, \{19,17,11,1\}, \{19,11,15,1\}, \{19,10,15,1\},
\end{eqnarray*}
by the symmetry of $2O$ and $M^{(k)}$.

\medskip
\noindent
{\bf Category 9.} $n=4$ and $i_1=19$, $i_2,i_3,i_4\in\mathcal{Q}_1$. Then $\Lambda^{(k)}=\{1,9,15,19\}$, obviously.

\medskip
\noindent
{\bf 3.2. The restriction between $\Lambda^{(k_1)}$ and $\Lambda^{(k_2)}$, for $k_1\neq k_2$}

By routine computation, we list some restricting conditions between different faces as follow:

\noindent
For an arbitrary point ${\bf x}\in Q^{(k)}_{1}$, we have
$$(2O+{\bf x})\supset \big(Q^{(k+1)}_{1}\cup Q^{(k+3)}_{1}\big).$$
Therefore, if $1\in\Lambda^{(k)}$, we must have
$$1\notin\Lambda^{(k+1)},\ 1\notin\Lambda^{(k+3)}.\eqno(30)$$

For an arbitrary point ${\bf x}\in Q^{(k)}_{2}$, we have
$$(2O+{\bf x})\supset \big(Q^{(k+3)}_{3}\cup Q^{(k+3)}_{4}\cup Q^{(k+3)}_{7}\cup Q^{(k+3)}_{8}\big).$$
Therefore, if $2\in\Lambda^{(k)}$, we must have
$$\Lambda^{(k+3)} \cap \{3,4,7,8\}= \emptyset.\eqno(31)$$

For an arbitrary point ${\bf x}\in Q^{(k)}_{4}$, we have
$$(2O+{\bf x})\supset \big(Q^{(k+1)}_{2}\cup Q^{(k+1)}_{4}\cup Q^{(k+1)}_{5}\cup Q^{(k+1)}_{6}\cup Q^{(k+3)}_{3}\cup Q^{(k+3)}_{4}\cup Q^{(k+3)}_{7}\cup Q^{(k+3)}_{8}\big).$$
Therefore, if $4\in\Lambda^{(k)}$, we must have
$$\Lambda^{(k+1)} \cap \{2,4,5,6\}= \emptyset,\ \Lambda^{(k+3)} \cap \{3,4,7,8\}= \emptyset.\eqno(32)$$

For an arbitrary point ${\bf x}\in \big(Q^{(k)}_{5}\cup Q^{(k)}_{6}\big)$, we have
$$(2O+{\bf x})\supset \big(Q^{(k+3)}_{3}\cup Q^{(k+3)}_{4}\cup Q^{(k+3)}_{7}\cup Q^{(k+3)}_{8}\cup Q^{(k+3)}_{16}\cup Q^{(k+3)}_{18}\cup Q^{(k+3)}_{19}\big).$$
Therefore, if $\Lambda^{(k)} \cap \{5,6\}\neq \emptyset$, we must have
$$\Lambda^{(k+3)} \cap \{3,4,7,8,16,18,19\}= \emptyset.\eqno(33)$$

For an arbitrary point ${\bf x}\in Q^{(k)}_{19}$, we have
$$(2O+{\bf x})\supset \big(Q^{(k+1)}_{5}\cup Q^{(k+1)}_{6}\cup Q^{(k')}_{13}\cup Q^{(k')}_{14}\cup Q^{(k+3)}_{7}\cup Q^{(k+3)}_{8}\big).$$
Therefore, if $19\in\Lambda^{(k)}$, we must have
$$\Lambda^{(k+1)} \cap \{5,6\}= \emptyset,\ \Lambda^{(k+2)} \cap \{13,14\}= \emptyset,\ \Lambda^{(k+3)} \cap \{7,8\}= \emptyset.\eqno(34)$$

By symmetry, we have:
\begin{align}
&{\rm if}\ 3\in\Lambda^{(k)},\ {\rm then}\ \Lambda^{(k+1)} \cap \{2,4,5,6\}= \emptyset,\tag{35}\\
&{\rm if}\ \Lambda^{(k)} \cap \{7,8\}\neq \emptyset,\ {\rm then}\ \Lambda^{(k+1)} \cap \{2,4,5,6,10,12,19\}= \emptyset,\tag{36}\\
&{\rm if}\ 9\in\Lambda^{(k)},\ {\rm then}\ 9\notin\Lambda^{(k+2)},\ 15\notin\Lambda^{(k+3)},\tag{37}\\
&{\rm if}\ 10\in\Lambda^{(k)},\ {\rm then}\ \Lambda^{(k+3)} \cap \{7,8,16,18\}= \emptyset,\tag{38}\\
&{\rm if}\ 11\in\Lambda^{(k)},\ {\rm then}\ \Lambda^{(k+2)} \cap \{11,12,13,14\}= \emptyset,\tag{39}\\
&{\rm if}\ 12\in\Lambda^{(k)},\ {\rm then}\ \Lambda^{(k+2)} \cap \{11,12,13,14\}= \emptyset,\ \Lambda^{(k+3)} \cap \{7,8,16,18\}= \emptyset,\tag{40}\\
&{\rm if}\ \Lambda^{(k)} \cap \{13,14\}\neq \emptyset,\ {\rm then}\ \Lambda^{(k+2)} \cap \{11,12,13,14,17,18,19\}= \emptyset,\tag{41}\\
&{\rm if}\ 15\in\Lambda^{(k)},\ {\rm then}\ 9\notin\Lambda^{(k+1)},\ 15\notin\Lambda^{(k+2)},\tag{42}\\
&{\rm if}\ 16\in\Lambda^{(k)},\ {\rm then}\ \Lambda^{(k+1)} \cap \{5,6,10,12\}= \emptyset,\tag{43}\\
&{\rm if}\ 17\in\Lambda^{(k)},\ {\rm then}\ \Lambda^{(k+2)} \cap \{13,14,17,18\}= \emptyset,\tag{44}\\
&{\rm if}\ 18\in\Lambda^{(k)},\ {\rm then}\ \Lambda^{(k+1)} \cap \{5,6,10,12\}= \emptyset,\ \Lambda^{(k+2)} \cap \{13,14,17,18\}= \emptyset.\tag{45}
\end{align}

For an arbitrary point ${\bf x}\in \big(Q^{(1)}_{4}\cup Q^{(1)}_{6}\cup Q^{(1)}_{19}\big)$, we have
$$(2O+{\bf x})\supset {\rm conv}\big\{(\mbox{$1,\frac{1}{3},\frac{2}{3}),(\frac{4}{3},0,\frac{2}{3}),(\frac{4}{3},\frac{1}{3},\frac{1}{3}),(1,\frac{2}{3},1),(\frac{4}{3},\frac{1}{3},1),(\frac{4}{3},\frac{2}{3},\frac{2}{3}),(\frac{5}{3},\frac{1}{3},\frac{2}{3}$})\big\}.$$
For an arbitrary point
$${\bf y}\in \Big(Q^{(1)}_{16}\setminus{\rm conv}\big\{(\mbox{$1,\frac{1}{3},\frac{2}{3}),(\frac{4}{3},0,\frac{2}{3}),(\frac{4}{3},\frac{1}{3},\frac{1}{3}),(1,\frac{2}{3},1),(\frac{4}{3},\frac{1}{3},1),(\frac{4}{3},\frac{2}{3},\frac{2}{3}),(\frac{5}{3},\frac{1}{3},\frac{2}{3}$})\big\}\Big),$$
we have
$$(2O+{\bf y})\supset \big(Q^{(2)}_{9}\cup Q^{(2)}_{11}\cup Q^{(1')}_{15}\big).$$
Therefore, if $\Lambda^{(1)} \cap \{4,6,19\}\neq \emptyset$ and $16\in\Lambda^{(1)}$, we must have
$$\Lambda^{(2)} \cap \{9,11\}= \emptyset,\ 15\notin\Lambda^{(3)}.$$

By symmetry, we have:
\begin{align}
&{\rm if}\ \Lambda^{(k)} \cap \{4,6,19\}\neq \emptyset\ {\rm and}\ 16\in\Lambda^{(k)},  \ {\rm then}\ \Lambda^{(k+1)} \cap \{9,11\}= \emptyset,\ 15\notin\Lambda^{(k+2)},\tag{46}\\
&{\rm if}\ \Lambda^{(k)} \cap \{6,12,19\}\neq \emptyset\ {\rm and}\ 17\in\Lambda^{(k)},  \ {\rm then}\ \Lambda^{(k+2)} \cap\{15,16\}= \emptyset,\ 9\notin\Lambda^{(k+1)},\tag{47}\\
&{\rm if}\ \Lambda^{(k)} \cap \{8,18,19\}\neq \emptyset\ {\rm and}\ 11\in\Lambda^{(k)} ,  \ {\rm then}\ \Lambda^{(k+2)} \cap \{9,10\}= \emptyset,\ 15\notin\Lambda^{(k+3)},\tag{48}\\
&{\rm if}\ \Lambda^{(k)} \cap \{12,14,19\}\neq \emptyset\ {\rm and}\ 2\in\Lambda^{(k)},  \ {\rm then}\ \Lambda^{(k+3)} \cap \{1,2\}= \emptyset,\ 1\notin\Lambda^{(k+1)},\tag{49}\\
&{\rm if}\ \Lambda^{(k)} \cap \{14,18,19\}\neq \emptyset\ {\rm and}\ 3\in\Lambda^{(k)},  \ {\rm then}\ \Lambda^{(k+1)} \cap \{1,3\}= \emptyset,\ 1\notin\Lambda^{(k+3)},\tag{50}\\
&{\rm if}\ \Lambda^{(k)} \cap \{4,8,19\}\neq \emptyset\ {\rm and}\ 10\in\Lambda^{(k)},  \ {\rm then}\ \Lambda^{(k+3)} \cap \{15,17\}= \emptyset,\ 9\notin\Lambda^{(k+2)}.\tag{51}
\end{align}

For an arbitrary point ${\bf x}\in Q^{(1)}_{13}$, we have
$$(2O+{\bf x})\supset {\rm conv}\big\{(\mbox{$\frac{1}{3},\frac{2}{3},1),(\frac{2}{3},\frac{1}{3},1),(\frac{2}{3},\frac{2}{3},\frac{2}{3}),(\frac{2}{3},\frac{2}{3},\frac{4}{3}),(\frac{2}{3},1,1),(1,\frac{2}{3},1)$}\big\}.$$
For an arbitrary point
$${\bf y}\in \Big(Q^{(1)}_{4}\setminus {\rm conv}\big\{(\mbox{$\frac{1}{3},\frac{2}{3},1),(\frac{2}{3},\frac{1}{3},1),(\frac{2}{3},\frac{2}{3},\frac{2}{3}),(\frac{2}{3},\frac{2}{3},\frac{4}{3}),(\frac{2}{3},1,1),(1,\frac{2}{3},1$})\big\}\Big),$$
we have
$$(2O+{\bf y})\supset (Q^{(2)}_{1}\cup Q^{(2)}_{3}\cup Q^{(4)}_{1}\cup Q^{(4)}_{2}).$$
Therefore, if $\{4,13\}\subset\Lambda^{(1)}$, we must have
$$\Lambda^{(2)} \cap\{1,3\}= \emptyset,\ \Lambda^{(4)} \cap\{1,2\}= \emptyset.$$

By symmetry, we have:
\begin{align}
&{\rm if}\ \{4,13\}\subset\Lambda^{(k)},  \ {\rm then}\ \Lambda^{(k+1)} \cap \{1,3\}= \emptyset,\ \Lambda^{(k+3)} \cap \{1,2\}= \emptyset,\tag{52}\\
&{\rm if}\ \{7,12\}\subset\Lambda^{(k)},  \ {\rm then}\ \Lambda^{(k+2)} \cap \{9,10\}= \emptyset,\ \Lambda^{(k+3)} \cap \{15,17\}= \emptyset,\tag{53}\\
&{\rm if}\ \{5,18\}\subset\Lambda^{(k)},  \ {\rm then}\ \Lambda^{(k+1)} \cap \{9,11\}= \emptyset,\ \Lambda^{(k+2)} \cap \{15,16\}= \emptyset.\tag{54}
\end{align}

For an arbitrary point ${\bf x}\in \big(Q^{(2)}_{1}\cup Q^{(4)}_{1}\big)$, we have
\begin{eqnarray*}
(2O+{\bf x})\supset {\rm conv}\hspace{-0.5cm} &\big\{\hspace{-0.5cm}&(\mbox{$0,\frac{1}{3},\frac{5}{3}),(\frac{1}{3},0,\frac{5}{3}),(0,\frac{2}{3},\frac{4}{3}),(\frac{2}{3},0,\frac{4}{3}),(\frac{1}{3},\frac{2}{3},1$}),\\
& &(\mbox{$\frac{2}{3},\frac{1}{3},1),(\frac{1}{3},\frac{2}{3},\frac{5}{3}),(\frac{2}{3},\frac{1}{3},\frac{5}{3}),(\frac{1}{3},1,\frac{4}{3}),(1,\frac{1}{3},\frac{4}{3}$})\big\}.
\end{eqnarray*}
For an arbitrary point ${\bf y}\in \big(Q^{(4)}_{15}\cup Q^{(1')}_{9}\big)$, we have
\begin{eqnarray*}
(2O+{\bf y}) \supset {\rm conv}\hspace{-0.5cm} &\big\{\hspace{-0.5cm}&(\mbox{$\frac{1}{3},\frac{5}{3},0),(0,\frac{5}{3},\frac{1}{3}),(\frac{2}{3},\frac{4}{3},0),(0,\frac{4}{3},\frac{2}{3}),(\frac{2}{3},1,\frac{1}{3}$}),\\
& &(\mbox{$\frac{1}{3},1,\frac{2}{3}),(\frac{2}{3},\frac{5}{3},\frac{1}{3}),(\frac{1}{3},\frac{5}{3},\frac{2}{3}),(1,\frac{4}{3},\frac{1}{3}),(\frac{1}{3},\frac{4}{3},1$})\big\}.
\end{eqnarray*}
For an arbitrary point
\begin{eqnarray*}
{\bf z}\in  \big(Q^{(1)}_{10}\cup  Q^{(1)}_{11}\cup  Q^{(1)}_{12}\big) \setminus {\rm conv}\hspace{-0.5cm} &\big\{\hspace{-0.5cm}&(\mbox{$\frac{1}{3},\frac{5}{3},0),(0,\frac{5}{3},\frac{1}{3}),(\frac{2}{3},\frac{4}{3},0),(0,\frac{4}{3},\frac{2}{3}),(\frac{2}{3},1,\frac{1}{3}$}),\\
& &(\mbox{$\frac{1}{3},1,\frac{2}{3}),(\frac{2}{3},\frac{5}{3},\frac{1}{3}),(\frac{1}{3},\frac{5}{3},\frac{2}{3}),(1,\frac{4}{3},\frac{1}{3}),(\frac{1}{3},\frac{4}{3},1$})\big\},
\end{eqnarray*}
we have
\begin{eqnarray*}
(2O+{\bf z})\supset\big(Q^{(1)}_{2}\cup Q^{(1)}_{3}\cup Q^{(1)}_{4}\big)\setminus{\rm conv}\hspace{-0.5cm} &\big\{\hspace{-0.5cm}&(\mbox{$0,\frac{1}{3},\frac{5}{3}),(\frac{1}{3},0,\frac{5}{3}),(0,\frac{2}{3},\frac{4}{3}),(\frac{2}{3},0,\frac{4}{3}),(\frac{1}{3},\frac{2}{3},1$}),\\
& &(\mbox{$\frac{2}{3},\frac{1}{3},1),(\frac{1}{3},\frac{2}{3},\frac{5}{3}),(\frac{2}{3},\frac{1}{3},\frac{5}{3}),(\frac{1}{3},1,\frac{4}{3}),(1,\frac{1}{3},\frac{4}{3}$})\big\}.
\end{eqnarray*}
Therefore, if $1\in\Lambda^{(2)}\ {\rm or}\ 1\in\Lambda^{(4)}$, and $9\in\Lambda^{(3)}\ {\rm or}\ 15\in\Lambda^{(4)}$, combining with (12)-(15), we have
$${\rm card}\big\{\Lambda^{(1)} \cap \{2,3,4,10,11,12\}\big\}\leq1.$$

By symmetry, we have:

\noindent
If $1\in\Lambda^{(k+1)}\ {\rm or}\ 1\in\Lambda^{(k+3)}$, and $9\in\Lambda^{(k+2)}\ {\rm or}\ 15\in\Lambda^{(k+3)}$, then
$${\rm card}\big\{\Lambda^{(k)} \cap \{2,3,4,10,11,12\}\big\}\leq1.\eqno(55)$$
If $1\in\Lambda^{(k+1)}\ {\rm or}\ 1\in\Lambda^{(k+3)}$, and $9\in\Lambda^{(k+1)}\ {\rm or}\ 15\in\Lambda^{(k+2)}$, then
$${\rm card}\big\{\Lambda^{(k)} \cap \{2,3,4,16,17,18\}\big\}\leq1.\eqno(56)$$
If $9\in\Lambda^{(k+1)}\ {\rm or}\ 15\in\Lambda^{(k+2)}$, and $9\in\Lambda^{(k+2)}\ {\rm or}\ 15\in\Lambda^{(k+3)}$, then
$${\rm card}\big\{\Lambda^{(k)} \cap \{10,11,12,16,17,18\}\big\}\leq1.\eqno(57)$$

For an arbitrary point ${\bf x}\in {\rm conv}\big\{(\frac{1}{3},\frac{1}{3},\frac{4}{3}),(\frac{1}{3},\frac{2}{3},1), (\frac{2}{3},\frac{1}{3},1),(\frac{2}{3},\frac{2}{3},\frac{4}{3})\big\},$ we have
$$(2O+{\bf x})\supset \big(Q^{(2)}_{1}\cup Q^{(2)}_{3}\cup Q^{(4)}_{1}\cup Q^{(4)}_{2}\big).$$

\noindent
For an arbitrary point ${\bf y}\in {\rm conv}\big\{(\frac{4}{3},\frac{1}{3},\frac{1}{3}),(\frac{5}{3},\frac{2}{3},\frac{1}{3}), (\frac{4}{3},\frac{2}{3},0),(\frac{5}{3},\frac{1}{3},0)\big\},$ we have
$$(2O+{\bf y})\supset \big(Q^{(2)}_{9}\cup Q^{(1')}_{15}\cup Q^{(1')}_{16}\big).$$

\noindent
For an arbitrary point ${\bf z}\in \Big(Q^{(1)}_{4}\setminus{\rm conv}\big\{(\frac{1}{3},\frac{1}{3},\frac{4}{3}),(\frac{1}{3},\frac{2}{3},1), (\frac{2}{3},\frac{1}{3},1),(\frac{2}{3},\frac{2}{3},\frac{4}{3})\big\}\Big),$ we have
$$(2O+{\bf z})\supset \Big(Q^{(1)}_{17}\setminus{\rm conv}\big\{(\mbox{$\frac{4}{3},\frac{1}{3},\frac{1}{3}),(\frac{5}{3},\frac{2}{3},\frac{1}{3}), (\frac{4}{3},\frac{2}{3},0),(\frac{5}{3},\frac{1}{3},0)$}\big\}\Big).$$

\noindent
Therefore, if $\{4,17\}\subset\Lambda^{(1)}$, we have
$$\Lambda^{(2)}\cap\{1,3\}=\emptyset,\ \Lambda^{(4)}\cap\{1,2\}=\emptyset,$$
or
$$9\notin\Lambda^{(2)},\ \Lambda^{(3)}\cap \{15,16\}=\emptyset.\eqno(58)$$

By symmetry, we have:

\noindent
If $\{4,11\}\subset\Lambda^{(1)}$, we have
$$\Lambda^{(2)}\cap\{1,3\}=\emptyset,\ \Lambda^{(4)}\cap\{1,2\}=\emptyset,$$
or
$$15\notin\Lambda^{(4)},\ \Lambda^{(3)}\cap\{9,10\}=\emptyset.\eqno(59)$$

Now we show that: a combination of $\Lambda^{(k)},\ k=1,2,3,4$ as where categorized before, cannot satisfy conditions (30)-(59).

\medskip
\noindent
{\bf Case 1.} Category 3 is used. Without loss of generality, we suppose
$$\Lambda^{(1)}=\{4,13\}.\eqno(1.1)$$

Since $\Lambda\cap P^{(2)}_1\neq\emptyset$ and $\Lambda\cap P^{(4)}_1\neq\emptyset,$ by (5), combining (1.1) with (6), (32) and (52), we have
$$\Lambda^{(2)} \cap \{7,8\}\neq \emptyset\eqno(1.2)$$
and
$$\Lambda^{(4)} \cap \{5,6\}\neq \emptyset.\eqno(1.3)$$

Combining (1.1), (1.2), (1.3) with (7), (33), (36) and (41), we have
$$\Lambda\cap P^{(3)}_2=\emptyset,$$
which contradicts (5). Therefore, Category 3 cannot be used, and
$$\Lambda^{(k)} \cap \{5,7,13\}=\emptyset$$
holds for all $k$.

\medskip
\noindent
{\bf Case 2.} Category 5 is used. Without loss of generality, we suppose
$$\Lambda^{(1)}\supset\{4,16\},\ \Lambda^{(1)} \cap \{9,10,11\}\neq \emptyset.\eqno(2.1)$$

Since $\Lambda\cap P^{(2)}_4\neq\emptyset$ by (5), combining (2.1) with (9), (43), (46) and the conclusion of Case 1, we  have
$$14\in\Lambda^{(2)}\eqno(2.2)$$
and
$$\Lambda^{(2)}\cap\{1,2,3\}\neq\emptyset,\eqno(2.3)$$
by the categorization before.

Since $\Lambda\cap P^{(4)}_5\neq\emptyset$ by (5), combining (2.1), (2.2) with (10), (32) and (41), we have
$$6\in\Lambda^{(4)}\eqno(2.4)$$
and
$$\Lambda^{(4)}\cap \{15,16\}\neq\emptyset,\eqno(2.5)$$
by the categorization before.

Suppose
$$16\in\Lambda^{(4)},\eqno(2.5.1)$$
combining with (2.4), (43) and (46), we have
$$\Lambda^{(1)} \cap \{9,10,11\}=\emptyset,$$
which contradicts (2.1). Therefore, we have
$$15\in\Lambda^{(4)},\eqno(2.5.2)$$
combining with (2.1) and (42), we have
$$\Lambda^{(1)} \cap \{10,11\}\neq \emptyset.\eqno(2.6)$$

If
$$1\in\Lambda^{(2)},\eqno(2.3.1)$$
combining with (2.5.2) and (55), we have
$${\rm card}\big\{\Lambda^{(1)} \cap \{4,10,11\}\big\}\leq1,$$
which contradicts (2.1) and (2.6).

If
$$2\in\Lambda^{(2)},\eqno(2.3.2)$$
by (31), we have $4\notin\Lambda^{(1)}$, which contradicts (2.1). Therefore, we have
$$3\in\Lambda^{(2)}.\eqno(2.3.3)$$

By (2.3.3), (2.2), (6), (35), (50) and the conclusion of Case 1, we have
$$8\in\Lambda^{(3)},$$
which contradicts (2.4) and (33). Therefore, Category 5 cannot be used.

\medskip
\noindent
{\bf Case 3.} Category 6 is used. Without loss of generality, we suppose
$$\Lambda^{(1)}\supset\{4,17\},\ \Lambda^{(1)} \cap \{9,11\}\neq \emptyset.\eqno(3.1)$$

Suppose
$$\Lambda^{(2)}\cap\{1,3\}=\emptyset,\ \Lambda^{(4)}\cap\{1,2\}=\emptyset.\eqno(3.2.1)$$
Since $\Lambda\cap P^{(2)}_1\neq\emptyset$ and $\Lambda\cap P^{(4)}_1\neq\emptyset$ by (5), combining with (3.1), (6), (32) and the conclusion of Case 1, we have
$$8\in\Lambda^{(2)}\eqno(3.2.1.1)$$
and
$$6\in\Lambda^{(4)}.\eqno(3.2.1.2)$$

Since $\Lambda\cap P^{(3)}_2\neq\emptyset$ by (5), combining with (7), (3.2.1.1), (3.2.1.2), (33) and (36), we have
$$14\in\Lambda^{(3)}.\eqno(3.2.1.3)$$
By (41), we have $17\notin\Lambda^{(1)}$, which contradicts (3.1).

Suppose (3.2.1) is not hold, then by (3.1) and (58), we have
$$9\notin\Lambda^{(2)},\ \Lambda^{(3)}\cap \{15,16\}=\emptyset.\eqno(3.2.2)$$
Since $\Lambda\cap P^{(3)}_{6}\neq\emptyset$ by (5), combining with (3.1), (11), (44) and the conclusion of Case 1, we have
$$8\in\Lambda^{(3)}\eqno(3.2.2.1)$$
and
$$\Lambda^{(3)}\cap\{9,10,11\}\neq\emptyset,\eqno(3.2.2.2)$$
by the categorization before.

By (3.2.2.1), (3.2.2.2), (37), (48) and (51), we have
$$9\notin\Lambda^{(1)},\eqno(3.2.2.3)$$
combining with (3.1), we have
$$11\in\Lambda^{(1)}.\eqno(3.2.2.4)$$

By (3.1), (3.2.2.4), (39) and (59), we have
$$\Lambda^{(3)}\cap\{9,10,11\}=\emptyset,$$
which contradicts (3.2.2.2). Therefore, Category 6 cannot be used.

\medskip
\noindent
{\bf Case 4.} Category 2 is used. Without loss of generality, we suppose
$$\Lambda^{(1)}=\{6,16\}.\eqno(4.1)$$

Since $\Lambda\cap P^{(2)}_{4}\neq\emptyset$ by (5), combining with (4.1), (9), (43), (46) and the conclusion of Case 1, we have
$$14\in\Lambda^{(2)}\eqno(4.2)$$
and
$$\Lambda^{(2)} \cap \{1,2,3\}\neq \emptyset,\eqno(4.3)$$
by the categorization before.

By (4.2), (4.3), (30), (49) and (50), we have
$$1\notin\Lambda^{(3)}.\eqno(4.4)$$
Since $\Lambda\cap P^{(4)}_{3}\neq\emptyset$ by (5), combining with (4.1), (4.2), (8), (33) and (41), we have
$$6\in\Lambda^{(4)}.\eqno(4.5)$$

Since $\Lambda\cap P^{(3)}_{1}\neq\emptyset$ by (5), combining with (4.4), (4.5), (6), (33) and the conclusion of Case 1, we have
$$\Lambda^{(3)}\cap \{2,6\}\neq\emptyset.\eqno(4.6)$$

If
$$2\in\Lambda^{(3)},\eqno(4.6.1)$$
since $\Lambda\cap P^{(3)}_{3}\neq\emptyset$ by (5), combining with (4.5), (8), (12) and (33), we have
$$\Lambda^{(3)}\cap\{12,14\}\neq\emptyset.$$
By the conclusion of Case 2, we have $12\notin\Lambda^{(3)}.$ Therefore, we have
$$14\in\Lambda^{(3)}.$$
Combining with (4.6.1), (31) and (49), we have
$$\Lambda^{(2)} \cap \{1,2,3\}=\emptyset,$$
which contradicts (4.3).

If
$$6\in\Lambda^{(3)},\eqno(4.6.2)$$
then we have
$$\Lambda^{(3)}\cap\{15,16,17\}\neq\emptyset,\eqno(4.6.2.1)$$
by the categorization before.
Combining with (4.1), (4.5), (33) and (46), we have
$$17\in\Lambda^{(3)},\eqno(4.6.2.2)$$
combining with (4.6.2) and (47), we have
$$16\notin\Lambda^{(1)},$$
which contradicts (4.1). Therefore, Category 2 cannot be used.

\medskip
\noindent
{\bf Case 5.} Category 8 is used. Without loss of generality, we suppose
$$\Lambda^{(1)}\supset\{16,19\},\ \Lambda^{(1)} \cap \{1,2,3\}\neq\emptyset,\ \Lambda^{(1)} \cap \{9,10,11\}\neq\emptyset.\eqno(5.1)$$

Since $\Lambda\cap P^{(2)}_{4}\neq\emptyset,$ combining with (5.1), (9), (43), (46) and the conclusion of Case 1, we have
$$14\in\Lambda^{(2)}\eqno(5.2)$$
and
$$\Lambda^{(2)} \cap\{1,2,3\}\neq\emptyset,\eqno(5.3)$$
by the categorization before.

By (5.2), (5.3), (30), (49) and (50), we have
$$1\notin\Lambda^{(1)}.\eqno(5.4)$$
Suppose
$$3\in\Lambda^{(1)},\eqno(5.1.1)$$
combining with (5.1), (35) and (50), we have
$$\Lambda^{(2)} \cap\{1,2,3\}=\emptyset,$$
which contradicts (5.3). Therefore, we have
$$3\notin\Lambda^{(1)},\eqno(5.1.2)$$
combining with (5.1) and (5.4), we have
$$2\in\Lambda^{(1)}.\eqno(5.5)$$

Since $\Lambda\cap P^{(4)}_{1}\neq\emptyset$ by (5), combining with (5.1), (5.5), (6), (31), (49) and the conclusion of Case 1, we have
$$6\in\Lambda^{(4)}\eqno(5.6)$$
and
$$\Lambda^{(4)} \cap \{15,16,17\}\neq\emptyset,\eqno(5.7)$$
by the categorization before.
By (5.1), (37), (48) and (51), we have
$$15\notin\Lambda^{(4)}.\eqno(5.8)$$

If
$$16\in\Lambda^{(4)},\eqno(5.9)$$
combining with (5.6), (43) and (46), we have
$$\Lambda^{(1)} \cap\{9,10,11\}=\emptyset,$$
which contradicts (5.1). Therefore $16\notin\Lambda^{(4)}$. Combining with (5.7) and (5.8), we have
$$17\in\Lambda^{(4)}.\eqno(5.10)$$
 
By (5.10) and (44), we have
$$14\notin\Lambda^{(2)},$$
which contradicts (5.2). Therefore, Category 8 cannot be used.

\medskip
\noindent
{\bf Case 6.} Category 9 is used. Without loss of generality, we suppose
$$\Lambda^{(1)}=\{1,9,15,19\}.\eqno(6.1)$$

By (6.1), (30), (34) and (42), we have
$$\Lambda^{(2)} \cap \{1,5,6,9\}=\emptyset.\eqno(6.2)$$
Therefore, $\Lambda^{(2)}$ cannot use Category 1, 4 and 9, $\Lambda^{(2)}$ must use Category 7. In this case,
$${\rm card}\big\{\Lambda^{(2)} \cap\{2,3,10,11\}\big\}\geq2,$$
which contradicts (6.1) and (55). Therefore, Category 9 cannot be used.

\medskip
\noindent
{\bf Case 7.} Category 7 is used. Without loss of generality, we suppose
$$\Lambda^{(1)}=\{2,11,16\}.\eqno(7.1)$$

Since $\Lambda\cap P^{(2)}_{4}\neq\emptyset$ by (5), combining with (7.1), (9), (43) and the conclusion of Case 1, we have
$$\Lambda^{(2)}\cap \{9,11,14\}\neq\emptyset.\eqno(7.2)$$

\medskip
\noindent
{\bf Case 7.1.}
$$9\in\Lambda^{(2)}.\eqno(7.2.1)$$

Since only Category 1, 4 and 7 is still available, we have only three options for $F^{(2)}$:

\medskip
\noindent
{\bf Case 7.1.1.}
If
$$\Lambda^{(2)}=\{8,9\},\eqno(7.2.1.1)$$
since $\Lambda\cap P^{(3)}_{4}\neq\emptyset$ by (5), combining with (7.1), (9), (36) and (39), we have
$$9\in\Lambda^{(3)}.$$
Combining with (7.2.1) and (57), we have
$${\rm card}\big\{\Lambda^{(1)} \cap\{10,11,12,16,17,18\}\big\}\leq1,$$
which contradicts (7.1).

\medskip
\noindent
{\bf Case 7.1.2.}
If
$$\Lambda^{(2)}=\{1,9,18\},\eqno(7.2.1.2)$$
by (56), we have
$${\rm card}\{\Lambda^{(1)} \cap\{2,3,4,16,17,18\}\}\leq1,$$
which contradicts (7.1).

\medskip
\noindent
{\bf Case 7.1.3.}
If
$$\Lambda^{(2)}=\{4,9,15\},\eqno(7.2.1.3)$$
since $\Lambda\cap P^{(3)}_{4}\neq\emptyset$ by (5), combining with (9), (7.1), (32), (39) and (42), we have
$$10\in\Lambda^{(3)}.$$
Since only Category 1, 4 and 7 is still available, by the categorization before, we have
$$\Lambda^{(3)}=\{3,10,17\}.$$

Since $\Lambda\cap P^{(4)}_{1}\neq\emptyset$ by (5), combining with (6), (7.1), (31) and (35), we have
$$1\in\Lambda^{(4)}.$$
Combining with (7.2.1.3) and (56), we have
$${\rm card}\{\Lambda^{(1)} \cap\{2,3,4,16,17,18\}\}\leq1,$$
which contradicts (7.1).

\medskip
\noindent
{\bf Case 7.2.}
$$11\in\Lambda^{(2)}.\eqno(7.2.2)$$

Since only Category 1, 4 and 7 is still available, by the categorization before, we have
$$\Lambda^{(2)}=\{2,11,16\}.\eqno(7.2.2.1)$$
Since $\Lambda\cap P^{(3)}_{4}\neq\emptyset$ by (5), combining with (7.1), (9), (39) and (43), we have
$$9\in\Lambda^{(3)},$$
which is same to Case 7.1, up to symmetry.

\medskip
\noindent
{\bf Case 7.3.}
$$14\in\Lambda^{(2)}.\eqno(7.2.3)$$ 

In this case, we have
$$\Lambda^{(2)}=\{1,14\},\eqno(7.2.3.1)$$
by the categorization before. Since $\Lambda\cap P^{(4)}_{5}\neq\emptyset$ by (5), combining with (7.1), (7.2.3), (10), (31) and (41), we have
$$6\in\Lambda^{(4)}.$$
Therefore, we have
$$\Lambda^{(4)}=\{6,15\},\eqno(7.2.3.2)$$
by the categorization before. By (7.2.3.1), (7.2.3.2) and (55), we have
$${\rm card}\{\Lambda^{(1)} \cap\{2,3,4,10,11,12\}\}\leq1,$$
which contradicts (7.1).

As a conclusion, Category 7 cannot be used.

\medskip
\noindent
{\bf Case 8.} Category 4 is used. Without loss of generality, we suppose
$$\Lambda^{(1)}=\{1,9,18\}.\eqno(8.1)$$

Since only Category 1 and 4 is still available, combining with (8.1), (30) and (37), we have
$$\Lambda^{(4)}=\{8,9\}.\eqno(8.2)$$
By (8.1), (8.2) and (37), we have
$$\Lambda^{(3)}=\{1,14\}.\eqno(8.3)$$
By (8.3) and (41), we have
$$18\notin\Lambda^{(1)},$$
which contradicts (8.1). Therefore, Category 4 cannot be used.

\medskip
\noindent
{\bf Case 9.} Category 1 used by all the faces. Without loss of generality, we suppose
$$\Lambda^{(1)}=\{1,14\}.\eqno(9.1)$$

Since $\Lambda\cap P^{(3)}_{5}\neq\emptyset$ by (5), combining with (9.1), (10) and (41), we have
$$\Lambda^{(3)}\cap\{6,8\}\neq\emptyset.\eqno(9.2)$$
Without loss of generality, we suppose
$$\Lambda^{(3)}=\{6,15\}.\eqno(9.3)$$

By (9.1), (9.3), (30) and (33), we have
$$\Lambda^{(2)}\cap\{1,8\}=\emptyset.$$
Therefore, we have
$$\Lambda^{(2)}=\{6,15\},\eqno(9.4)$$
by the categorization before. By (9.1), (9.3), (9.4), (30) and (42), we have
$$\Lambda^{(4)}\cap \{1,9,15\}=\emptyset,$$
which is a contradiction, since Category 1 cannot be used in $\Lambda^{(4)}$.

As a conclusion, for lattice packing $O+\Lambda$, (5) cannot hold, which means $\gamma(O,\Lambda)\geq\frac{7}{6}$ holds for all lattice sets. Particularly, since the lattice $\Lambda$ generated by ${\bf a}_1=(-\frac{2}{3},1,\frac{1}{3})$, ${\bf a}_2=(\frac{1}{3},-\frac{2}{3},1)$, ${\bf a}_3=(1,\frac{1}{3},-\frac{2}{3})$ given in \cite{Mi01} (see also \cite{Be01}) can be easily verified to satisfy that $\gamma(O,\Lambda)=\frac{7}{6}$, thus we have
$$\gamma^{*}(O)=7/6,$$
Theorem 1 is proved.
\end{proof}

\vspace{0.8cm}
\noindent
{\LARGE\bf 4. Several Examples about Problem 1}

\bigskip
In $\mathbb{E}^2$, it is known that the density of the densest lattice packing of a smoothed octagon is smaller than the density of the densest packing of a circular disk (see \cite{Ha01}) and the simultaneous packing-covering constant of a regular octagon is bigger than the simultaneous packing-covering constant of a circular disk (see \cite{Zo01,Zo03}).  However, in $\mathbb{E}^3$, evidences support Ulam's conjecture (see \cite{Co01}) which claims that the density of the densest packing of a convex body attains its minimum at balls. In this section, we present some examples about the simultaneous packing-covering analogy of Ulam's conjecture.

Suppose $C$ is a centrally symmetric convex body in $\mathbb{E}^3$. Let $\Lambda$ be a lattice generated by $\{{\bf a}_1, {\bf a}_2, {\bf a}_3\}$, let $V$ denote the set $\{{\bf o}, {\bf a}_1, {\bf a}_2, {\bf a}_1+{\bf a}_2, {\bf a}_3, {\bf a}_3-{\bf a}_2\}$, and let $P$ denote the convex hull of $V$. Then, we have the following criterion:

\medskip
\noindent
{\bf Lemma 4 (Zong \cite{Zo02}).} If $P\subset C+V,$ then $C+\Lambda$ is a lattice covering.

\bigskip
In 2000, Betke and Henk \cite{Be01} discovered an algorithm by which one can determine the density of the densest lattice packing of any given three-dimensional convex polytope. As an application they calculated densest lattice packings of all regular and Archimedean polytopes. Applying this criterion to Betke and Henk's discoveries, we obtain the following examples:

\medskip
\noindent
{\bf Example 1.} Let $C_1$ be the octahedron defined by $\{ (x_1,x_2,x_3):\ |x_1|+|x_2|+|x_3|\le 1\}$ and let $\Lambda_1$ denote the lattice with a basis ${\bf a}_1=(2/3,1,1/3)$, ${\bf a}_2=(-1/3,-2/3,1)$ and ${\bf a}_3=(-1,1/3,-2/3)$. One can prove that $C_1+\Lambda_1$ is a packing in $\mathbb{E}^3$ and ${7\over 6}C_1+\Lambda_1$ is a covering of $\mathbb{E}^3$. Therefore, we have
$$\gamma^{*}(C_1)\le \frac{7}{6}<\sqrt{5/3}=\gamma^*(B^3).$$

\medskip
\noindent
{\bf Example 2.} We take $\tau=\frac{\sqrt5+1}{2}$ and define
$$C_2=\{(x_1,x_2,x_3): |\tau x_1|+|x_2|\leq1,\ |\tau x_2|+|x_3|\leq1,\ |\tau x_3|+|x_1|\leq1\},$$
\begin{eqnarray*}
C_5=\{(x_1,x_2,x_3)\hspace{-0.3cm}&:&\hspace{-0.2cm} |x_1|+|x_2|+|x_3|\leq1,\ |\tau x_1|+|\mbox{$\frac{1}{\tau}$}x_3|\leq1,\ |\tau x_2|+|\mbox{$\frac{1}{\tau}$}x_1|\leq1,\\
& &\hspace{-0.2cm}|\tau x_3|+|\mbox{$\frac{1}{\tau}$}x_2|\leq1\},
\end{eqnarray*}
$$C_3=C_2\cap C_5,$$
and
$$C_4=C_2\cap\left(\mbox{$\frac{7+12\tau}{(3+4\tau)(1+\tau)}$}C_5\right).$$
Usually, $C_2$, $C_5$, $C_3$ and $C_4$ are called a dodecahedron, an icosahedron, an icosidodecahedron and a truncated dodecahedron. Let $\Lambda_2=\Lambda_3=\Lambda_4$ to be the lattice with a basis ${\bf a}_1=(0,\frac{2}{1+\tau},\frac{2}{1+\tau})$, ${\bf a}_2=(\frac{2}{1+\tau},0,\frac{2}{1+\tau})$ and ${\bf a}_3=(\frac{2}{1+\tau},\frac{2}{1+\tau},0).$ One can prove that, for $i=2$, $3$ and $4$, $C_i+\Lambda_i$ is a packing in $\mathbb{E}^3$ and $(\sqrt{5}-1)C_i+\Lambda_i$ is a covering of $\mathbb{E}^3$. Therefore, we have
$$\gamma^{*}(C_i)\le \sqrt{5}-1<\sqrt{5/3}=\gamma^*(B^3), \qquad i=2, 3, 4.$$

\medskip
\noindent
{\bf Example 3.} We continue to use the notation of Example 2 and define
$$C_6=C_5\cap\left(\mbox{$\frac{\frac{4}{3}+\tau}{1+\tau}$}C_2\right).$$
Usually, $C_6$ is called a truncated icosahedron. Let $\Lambda_5=\Lambda_6$ to be the lattice with a basis
${\bf a}_1=(\frac{4}{3},0,0)$, ${\bf a}_2=(0,\frac{4}{3},0)$ and ${\bf a}_3=(\frac{2}{3},\frac{2}{3},\frac{2}{3}).$
One can prove that, for $i=5$ and $6$, $C_i+\Lambda_i$ is a packing in $\mathbb{E}^3$ and $\sqrt{5/3}C_i+\Lambda_i$ is a covering of $\mathbb{E}^3$.
Therefore, we have
$$\gamma^{*}(C_i)\le \sqrt{5/3}= \gamma^*(B^3),\qquad i=5, 6.$$

\medskip
\noindent
{\bf Remark 1.} It is interesting to notice that, for an icosahedron, the optimal lattices for packing density are no longer optimal for simultaneous packing-covering constant. It is well-known (see \cite{Co01,Zo04}) that, for the unit ball $B^3$, the optimal lattices for packing density are different from the optimal lattices for covering density which are identical with the optimal lattices for simultaneous packing-covering constant.

\medskip
\noindent
{\bf Example 4.} Let $C_0$ denote the cube $\{(x_1,x_2,x_3):\ |x_1|,|x_2|,|x_3|\leq1\}.$ We define
$$C_7=C_0\cap (2C_1).$$
Usually $C_7$ is called a cubeoctahedron. Let $\Lambda_7$ denote the lattice with a basis ${\bf a}_1=(2,-\frac{1}{3},-\frac{1}{3})$, ${\bf a}_2=(-\frac{1}{3},2,-\frac{1}{3})$ and ${\bf a}_3=(-\frac{1}{3},-\frac{1}{3},2).$ One can prove that $C_7+\Lambda_7$ is a packing in $\mathbb{E}^3$ and ${7\over 6}C_7+\Lambda_7$ is a covering of $\mathbb{E}^3$. Therefore, we have
$$\gamma^{*}(C_7)\le {7\over 6}<\sqrt{5/3}= \gamma^*(B^3).$$

\medskip
\noindent
{\bf Example 5.} We define
\begin{eqnarray*}
C_8\hspace{-0.3cm}&=\hspace{-0.3cm}&\{(x_1,x_2,x_3):|x_1|+|x_2|\leq2+3\sqrt2,\ |x_1|+|x_3|\leq2+3\sqrt2,\ |x_2|+|x_3|\leq2+3\sqrt2\}\\
& &\cap(2\sqrt2+1)C_0\cap(3\sqrt2+3)C_1.
\end{eqnarray*}
Usually $C_8$ is called a truncated cubeoctahedron. Let $\Lambda_8$ denote the lattice with a basis ${\bf a}_1=(7.6568\cdots,\ -2.0339\cdots,\ 2.0339\cdots),$ ${\bf a}_2=(1.5185\cdots,\ 0.6901\cdots,\ 7.6568\cdots)$ and ${\bf a}_3=(6.1383\cdots,\ 5.6228\cdots,\ 2.7241\cdots).$ One can prove that $C_8+\Lambda_8$ is a packing in $\mathbb{E}^3$ and $\sqrt{5/3}C_8+\Lambda_8$ is a covering of $\mathbb{E}^3$. Therefore, we have
$$\gamma^{*}(C_8)\le \sqrt{5/3}=\gamma^*(B^3).$$

\medskip
\noindent
{\bf Example 6.} We define
$$C_9=\{(x_1,x_2,x_3):|x_1|+|x_2|\leq2,\ |x_1|+|x_3|\leq2,\ |x_2|+|x_3|\leq2\}\cap\sqrt2C_0\cap(4-\sqrt2)C_1.$$
Usually, it is called a rhombic cubeoctahedron. Let $\Lambda_9$ denote the lattice generated by ${\bf a}_1=(0,2,2)$, ${\bf a}_2=(2,0,2)$ and ${\bf a}_3=(2,2,0).$ It was shown by \cite{Be01} that the density of the densest lattice packing of $C_9$ is attained at $C_9+\Lambda_9$. However, the simultaneous packing-covering constant of $C_9+\Lambda_9$ is $\sqrt{2}$, which is much bigger than $\sqrt{5/3}$. On the other hand, let $\Lambda_9^*$ denote the body cubic center lattice generated by ${\bf a}_1=(\frac{4(4-\sqrt2)}{3},0,0)$, ${\bf a}_2=(0,\frac{4(4-\sqrt2)}{3},0)$ and ${\bf a}_3=(\frac{2(4-\sqrt2)}{3},\frac{2(4-\sqrt2)}{3},\frac{2(4-\sqrt2)}{3}),$ it can be verified that the simultaneous packing-covering constant of $C_9+\Lambda_9^*$ is between $\sqrt{5/3}$ and $\sqrt{2}$.

\vspace{0.6cm}\noindent
{\bf Acknowledgement.} This work is supported by the National Natural Science Foundation of China (NSFC12226006, NSFC11921001) and the Natural Key Research and Development Program of China (2018YFA0704701).

\vspace{0.8cm}\noindent
Yiming Li, Center for Applied Mathematics, Tianjin University, Tianjin 300072, P. R. China

\noindent
xiaozhuang@tju.edu.cn

\bigskip\noindent
Chuanming Zong, Center for Applied Mathematics, Tianjin University, Tianjin 300072, P. R. China

\noindent
cmzong@tju.edu.cn


\vspace{0.8cm}
\begin{thebibliography}{1}
{\large
\bibitem{Be01}U. Betke and M. Henk, Densest lattice packings of 3-polytopes, \emph{Comput. Geom.} \textbf{16} (2000), 157-186.
\bibitem{Bo01}K. B\"{o}r\"{o}czky, Closest packing and loosest covering of the space with balls, \emph{Studia Sci. Math. Hungar} \textbf{21} (1986), 79-89.
\bibitem{Pa01}P. Brass, W. Moser, J. Pach, \emph{Research Problems in Discrete Geometry}, Springer-Verlag, New York, 2005.
\bibitem{Bu01}G. L. Butler, Simultaneous packing and covering in Euclidean space, \emph{Proc. London Math. Soc.} \textbf{25} (1972), 721-735.
\bibitem{Co01}J. H. Conway and S. Torquato, Packing, Tiling, and Covering with Tetrahedra, {\it Proc. Natl. Acad.
Sci. USA} {\bf 103} (2006), 10612-10617.
\bibitem{Fe01}L. Fejes T\'{o}th, Close packing and loose covering with balls, \emph{Publ. Math. Debrecen} \textbf{23} (1976), 323-326.
\bibitem{Fe02}L. Fejes T\'{o}th, Remarks on the closest packing of convex discs, \emph{Comment. Math. Helv.} \textbf{53} (1978), 536-541.
\bibitem{Ha01}T. C. Hales, On the Reinhardt conjecture, {\it Vietnam J. Math.} {\bf 39} (2011), 287-307.
\bibitem{He01}M. Henk, \emph{Finite and Infinite Packings}, Habilitationsschrift, Universit\"at Siegen, 1995.
\bibitem{He02}M. Henk, A note on lattice packings via lattice refinements, {\it Expe. Math.} {\bf 27} (2018), 1-9.
\bibitem{Ho01}J. Horv\'{a}th, On close lattice packing of unit spheres in the space $E^n$, \emph{Proc. Steklov Math. Inst.} \textbf{152} (1982), 237-254.
\bibitem{Ho02}J. Horv\'{a}th, \emph{Several Problems of n-dimensional Discrete Geometry}, Ph.D. Thesis, Steklov Math. Inst., Moskow, 1986.
\bibitem{Li01}J. Linhart, Closest packings and closest coverings by translates of a convex disc, \emph{Studia Math. Hungar.} \textbf{13} (1978), 157-162.
\bibitem{Mi04}D. Micciancio, Almost perfect lattices, the covering radius problem, and applications to Ajtai's connection factor, {\it SIAM J. Comput.} {\bf 34} (2004), 118-169.
\bibitem{Mi01}H. Minkowski, Dichteste gitterf\"{o}rmige Lagerung kongruenter K\"{o}rper, \emph{Nachr. K. Ges. Wiss. G\"{o}ttingen} (1904), 311-355.
\bibitem{Ro01}C. A. Rogers, A note on coverings and packings, \emph{J. London Math. Soc.} \textbf{25} (1950), 327-331.
\bibitem{Ro02}C. A. Rogers, Lattice coverings of space: The Minkowski-Hlawka theorem, \emph{Proc. London Math. Soc.} (3)\textbf{8} (1958), 447-465.
\bibitem{Si01}C. L. Siegel, A mean value theorem in the geometry of numbers, \emph{Ann. of Math.} (2)\textbf{46} (1945), 340-347.
\bibitem{Zo01}C. Zong, Simultaneous packing and covering in the Euclidean plane, \emph{Monatsh. Math.} \textbf{134} (2002), 247-255.
\bibitem{Zo02}C. Zong, Simultaneous packing and covering in three-dimensional Euclidean space, \emph{J. London Math. Soc.} \textbf{67} (2003), 29-40.
\bibitem{Zo03}C. Zong, The simultaneous packing and covering constants in the plane, \emph{Adv. Math.} \textbf{218} (2008), 653-672.
\bibitem{Zo04}C. Zong, From deep holes to free planes, \emph{Bull. Amer. Math. Soc.} \textbf{39} (2002), 533-555.
\bibitem{Zo05}C. Zong, \emph{Sphere Packings}, Springer-Verlag, New York, 1999.}
\end{thebibliography}
\end{document}